\documentclass[journal,transmag]{IEEEtran}
\usepackage{cite}

\usepackage{color,soul}

\ifCLASSINFOpdf
\else
\fi
\usepackage{amsmath}
\usepackage{graphicx}
\usepackage{multirow}
\usepackage{tabu} 

\ifCLASSOPTIONcompsoc
  \usepackage[caption=false,font=normalsize,labelfont=sf,textfont=sf]{subfig}
\else
  \usepackage[caption=false,font=footnotesize]{subfig}
\fi
\begin{document}
%
\title{On Overcoming the Transverse Boundary Error of the SU/PG Scheme for Moving Conductor Problems\textsc{}}
\author{\IEEEauthorblockN{Sethupathy Subramanian\IEEEauthorrefmark{$\dag$}, 
Udaya Kumar,\IEEEauthorrefmark{$\ddag$} \textit{Senior Member, IEEE} and Sujata Bhowmick\IEEEauthorrefmark{$*$}}
\IEEEauthorblockA{\IEEEauthorrefmark{$\dag$}Department of Physics, University of Notre Dame, IN 46556, USA \\ 
\IEEEauthorrefmark{$\ddag$}Department of  Electrical Engineering, 
\IEEEauthorrefmark{$*$}Department of Electronic Systems Engineering, \\ Indian Institute of Science, Bangalore 560012, India }
\thanks{© 2021 IEEE. Personal use of this material is permitted.  Permission from IEEE must be obtained for all other uses, in any current or future media, including reprinting/republishing this material for advertising or promotional purposes, creating new collective works, for resale or redistribution to servers or lists, or reuse of any copyrighted component of this work in other works.}}


\markboth{IEEE Transactions on Magnetics$~~$-$~~$Accepted article $~~~~~~~~$ DOI 10.1109/TMAG.2021.3123748}%
{Shell \MakeLowercase{\textit{et al.}}: Bare Demo of IEEEtran.cls for IEEE Transactions on Magnetics Journals}
%
\IEEEtitleabstractindextext{%
\begin{abstract}
Conductor moving in magnetic field is quite common in electrical equipment. 
The numerical simulation of such problem is vital in their design and analysis 
of electrical equipment. The Galerkin finite element method (GFEM) is a commonly 
employed simulation tool, nonetheless, due to its inherent numerical instability 
at higher velocities, the GFEM requires upwinding techniques to handle moving 
conductor problems. The Streamline Upwinding/Petrov-Galerkin (SU/PG) scheme is a 
widely acknowledged upwinding technique, despite its error-peaking at the 
transverse boundary. This error at the transverse-boundary, is found to be 
leading to non-physical solutions. Several remedies have been suggested in the 
allied fluid dynamics literature, which employs non-linear, iterative techniques. 
The present work attempts to address this issue, by retaining the computational 
efficiency of the GFEM. By suitable analysis, it is shown that the source of the 
problem can be attributed to the Coulomb's gauge. Therefore, to solve the problem, 
the Coulomb's gauge is taken out from the formulation and the associated weak 
form is derived. The effectiveness of this technique is demonstrated with 
pertinent numerical results.
\end{abstract}

}

\maketitle

\IEEEdisplaynontitleabstractindextext

%
\IEEEpeerreviewmaketitle

\section{Introduction}
%
%
%
%
\IEEEPARstart{T}{he} moving conductor problems are the class of problems, 
{involving conductor moving in a magnetic field ($\bf B_a$). As a consequence, there would be an induced current ($\bf J$) in the conductor, which in turn produces a reaction magnetic field ($\bf b$).} The {quantification of the} reaction magnetic field and the associated induced currents, are vital for the design and analysis of the electrical equipment such as electromagnetic flowmeter, linear induction motor, eddy current brakes etc. For the ease of further discussion, the governing equations of the moving conductor problem  are enumerated below \cite{reducedb, fmbase},
\begin{equation} \label{eqge1}
\sigma \nabla \phi ~-~  (\nabla \cdot \dfrac{1}{\mu} \nabla) {\bf{A}} - \sigma~ {\bf{u}} \times \nabla \times {\bf{A}} = \sigma~ {\bf{u}} \times {\bf{B_{a}}}
\end{equation}	
\begin{equation} \label{eqge2}
\begin{split}
\nabla \cdot (\sigma \nabla \phi) - \nabla \cdot (\sigma ~ {\bf{u}} \times \nabla \times {\bf{A}}) =  \nabla \cdot (\sigma ~ \bf{u} \times \bf{B_{a}})
\end{split}
\end{equation}	
where, $\phi$ is the scalar potential arising out of the current flow, $\bf A$ is the magnetic vector potential associated with reaction magnetic field $\bf b$,  $\bf u$ is the velocity of the moving conductor, $\mu$  is the magnetic permeability and $\sigma$ is the electrical conductivity.

Owing to the coupled nature as well as the vectorial form of the governing equations, finding an analytical solution is nearly impossible and hence, numerical simulation is the only choice available.  Among the available numerical schemes, the Galerkin finite element method (GFEM) is widely employed for the simulation in different disciplines of engineering including electrical engineering. Despite its widespread success, as the velocities become higher the GFEM results {suffer from} numerical oscillations \cite{cdbook, upfdm1_sp, upfdm1_ru, up1}. Specifically, when the element Peclet number $Pe = \mu \sigma \|{\bf u}\| \Delta l/2 > 1$ ($\Delta l$ is the element length along the direction of the velocity), the numerical oscillations appears in the solution. This issue is comprehensively discussed in the allied literature and several remedies have been suggested  \cite{cdbook}.

The upwinding schemes are the earliest and popular ones suggested in the fluid dynamics literature and these have been extensively adopted in the moving conductor literature as well  \cite{mc6tf1, mc6tf2, mc2av1, mc3eb1, mc4ge1, mc5mc1,  mc2av2, mc3eb2, mcsupg_mcfit, mcsupg_cable, mcsupg_mfluid}. 
{It can be worth noting here that the occurrence of numerical oscillations at higher velocities is common to many numerical schemes including the edge element method and there are few upwinding remedies suggested in the literature} \cite{edgeup2, r:edgeup3}. {Owing to the presence of substantial cross-research, the upwinding schemes of nodal elements are well studied and well established. Hence, the present work employs nodal finite elements for the study.} Among the nodal upwinding schemes, the Streamline upwinding/Petrov Galerkin (SU/PG) scheme is the latest and the widely adopted {one}. Though, the SU/PG scheme successfully removes the numerical oscillations, there are numerous instances {wherein, it leads to} error at the transverse boundary \cite{discop}. This problem is scarcely discussed in the electrical engineering literature, {except for a recent effort in which the {existence} of such an error is pointed out for moving conductor problems} \cite{nemosu}. {Further, it is also shown that this} in turn leads to non-physical currents in the adjacent air region.

In order to {address this issue}, several approaches have been suggested in fluid dynamics literature \cite{discop1, discop, fic2}, which are termed as, 'Spurious Oscillations at Layers Diminishing (SOLD) methods'  \cite{soldreview1, soldreview2}. These schemes being non-linear involve an iterative calculation procedure, which is a computationally time consuming process. Moreover, the implementation of SOLD schemes for electrical engineering problems {can be rarely seen} in the literature. In any case, a recent work {has shown} that SOLD schemes are not fully free-from boundary error \cite{soldreview2}. 
{The present work aims to investigate the origin of this error and perhaps suggest a solution.}


Given these, the present work aims to investigate on the genesis of this boundary error in the SU/PG formulation for the moving conductor problem. It further aims to provide a simple and computationally efficient route to eliminate this error for moving conductor problems.


\section{Present Work}


\subsection{Analysis with 2D version of the problem}
\begin{figure}[!h]
\centering
\includegraphics[width=2.6in]{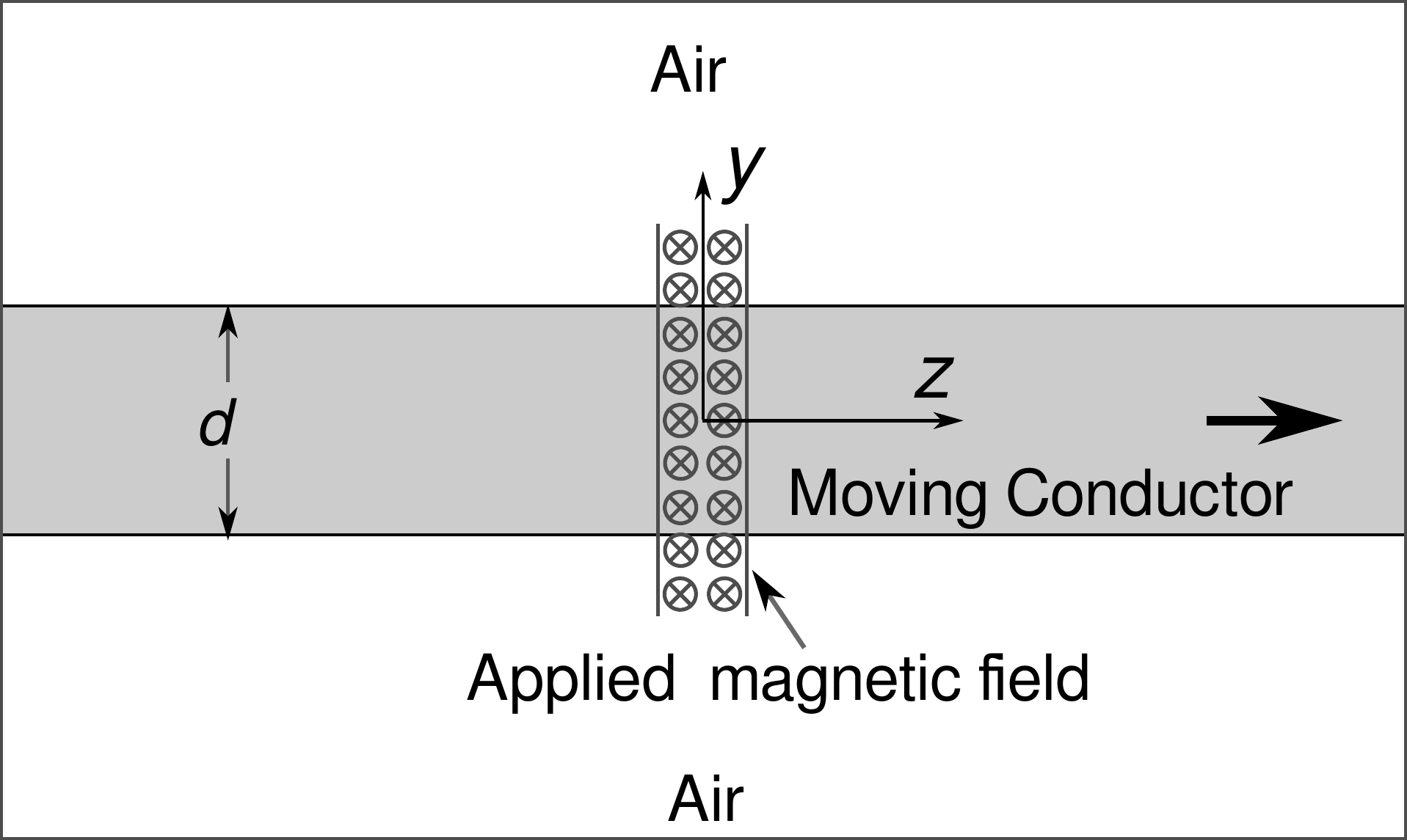}
\caption{Schematic of the 2D problem.}
\label{fig_sim}
\end{figure}

In reference \cite{nemosu}, 2D version of the moving conductor problem,
which has adequate representation of the 3D problem,
has been considered for evaluating the error at the interface. 
For the sake of completeness, a brief description of the 2D problem is 
provided here (see Fig.  \ref{fig_sim}). An infinite conducting slab with thickness $d$ is placed between two air layers of thickness $4d$. The conductor has a uniform velocity $u_z$ along the $z-$axis and a localized $x$-directed magnetic field is considered. Boundary in the flow direction is introduced at $20d$ on either side of the magnetic field. The material properties of the slab are $\sigma=7.21\times10^6Sm^{-1}$, $\mu_r=1$, $d=0.5m$. The FEM discretisation involved 5760 nodal, bilinear quadrilateral elements of rectangular shape. The element size is varied along the $y$-direction to accurately capture the currents circulating in the conductor. It can be noted that the 2D problem here represents the circulation of current and thus it consists of 2-components of vector potential $A_y, A_z$. 

The same 2D problem will be considered here for further analysis. In \cite{nemosu}, it is shown that only SU/PG scheme suffers from boundary error and {not the original Galerkin scheme or its variants} \cite{su1}. In order to trace the genesis of this error numerically, computed results are analysed to detect the presence of any inconsistency.

\begin{figure*}
		\centering
		\mbox{\subfloat[]{\label{divsA} 
\includegraphics[scale=0.44]{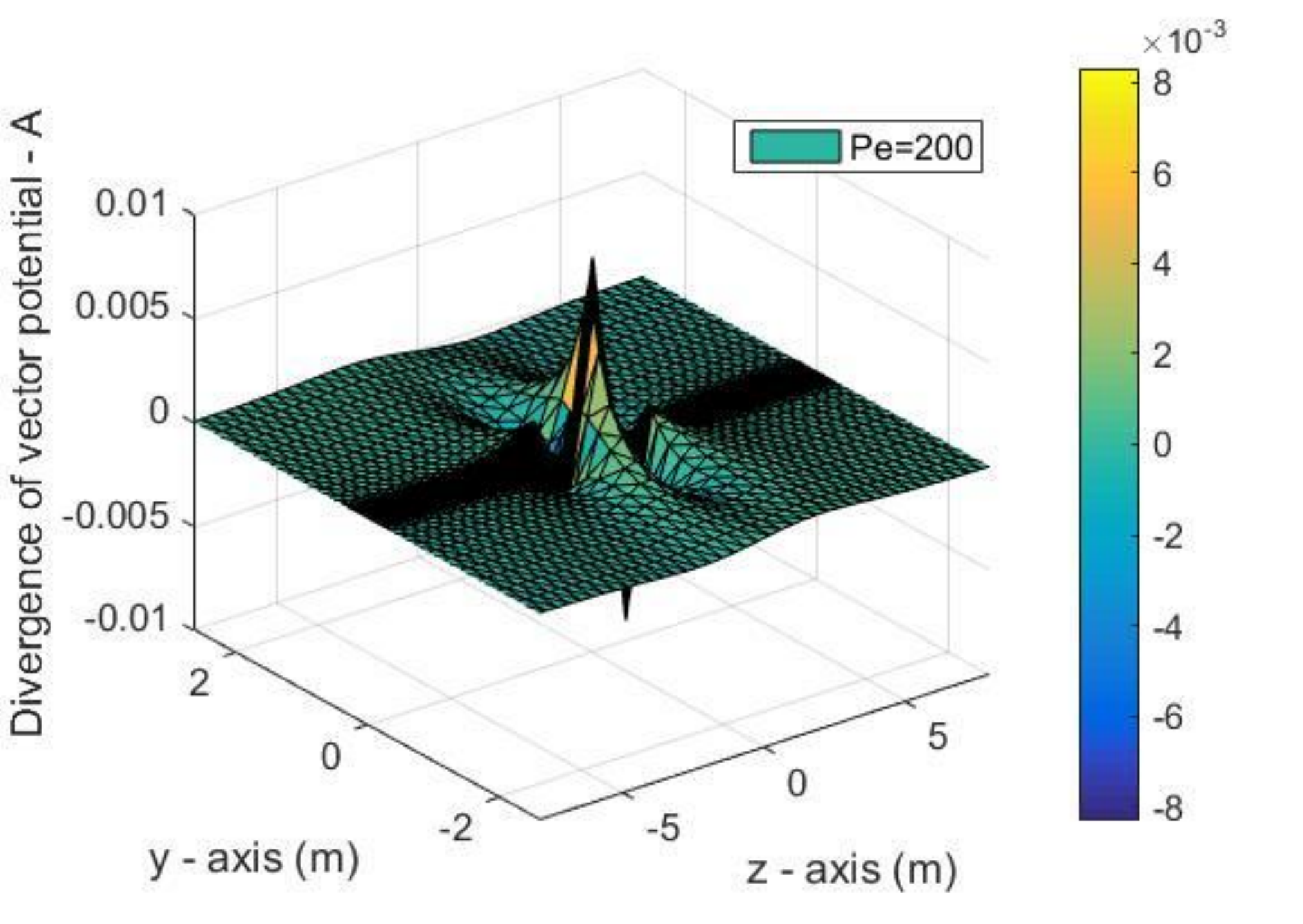}}}
\mbox{\subfloat[]{\label{divaA} 
\includegraphics[scale=0.44]{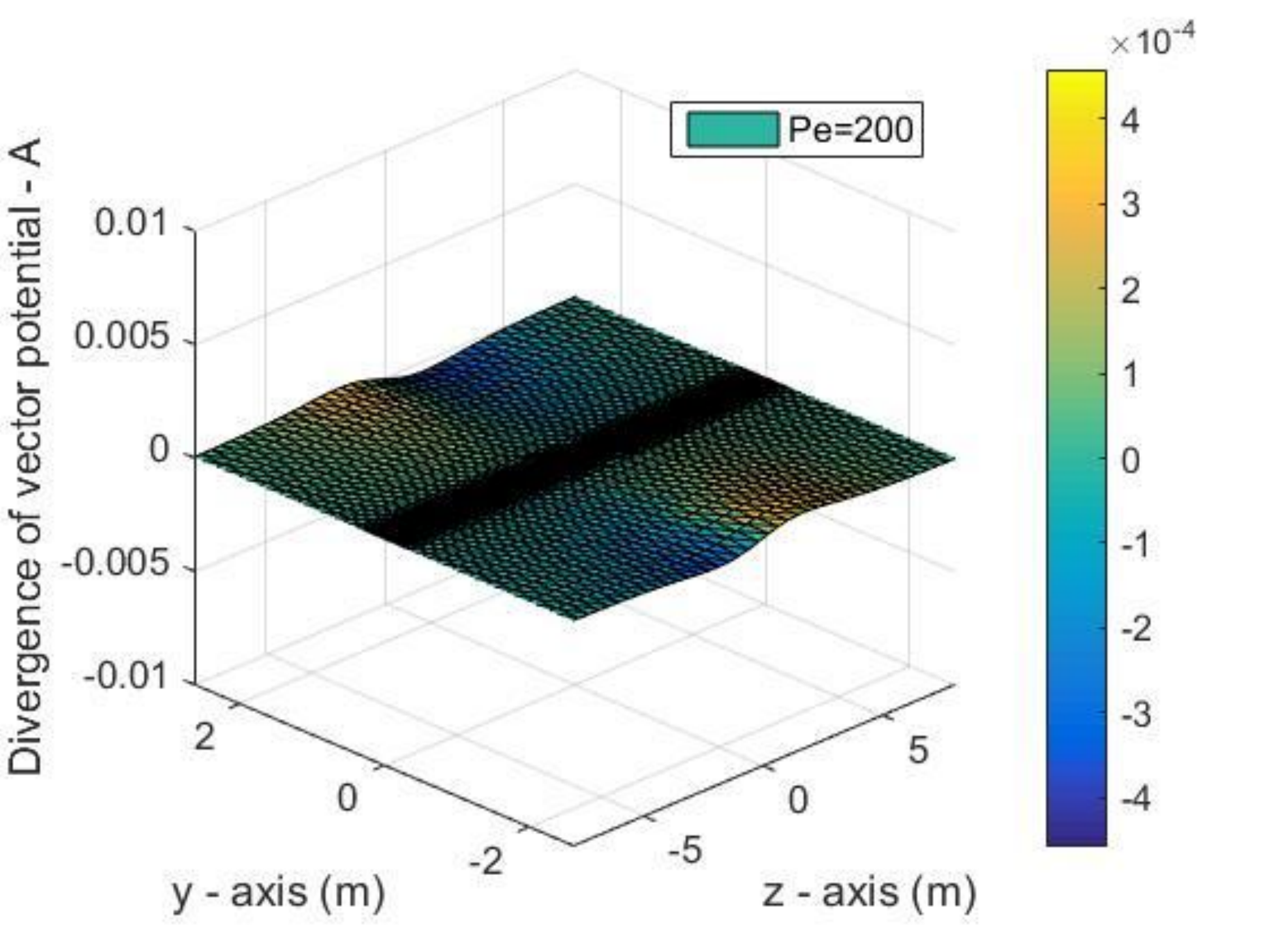}}}		
		\mbox{\subfloat[]{\label{dAerr} 
\includegraphics[scale=0.9]{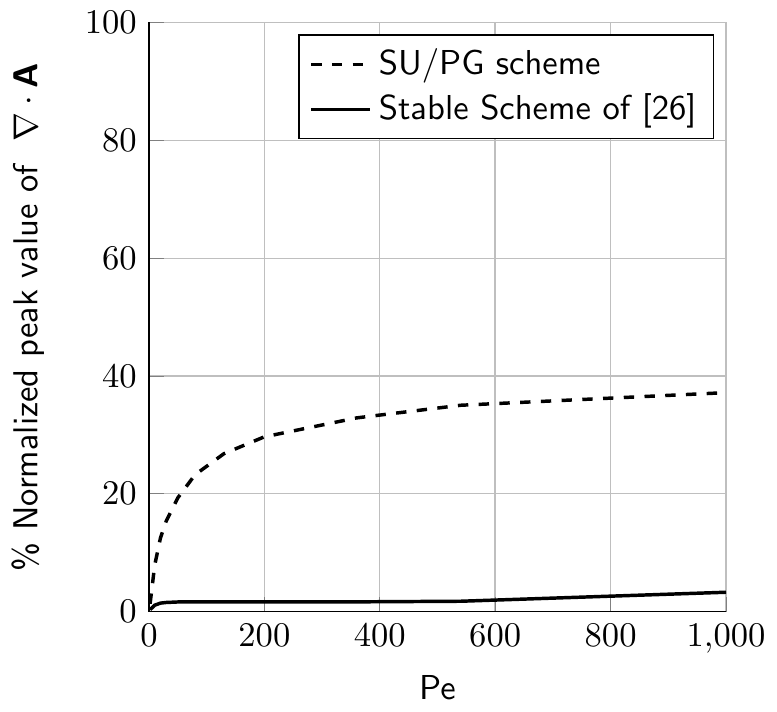}}} \hspace{1cm}
		\mbox{\subfloat[]{\label{Bderr} 
\includegraphics[scale=0.9]{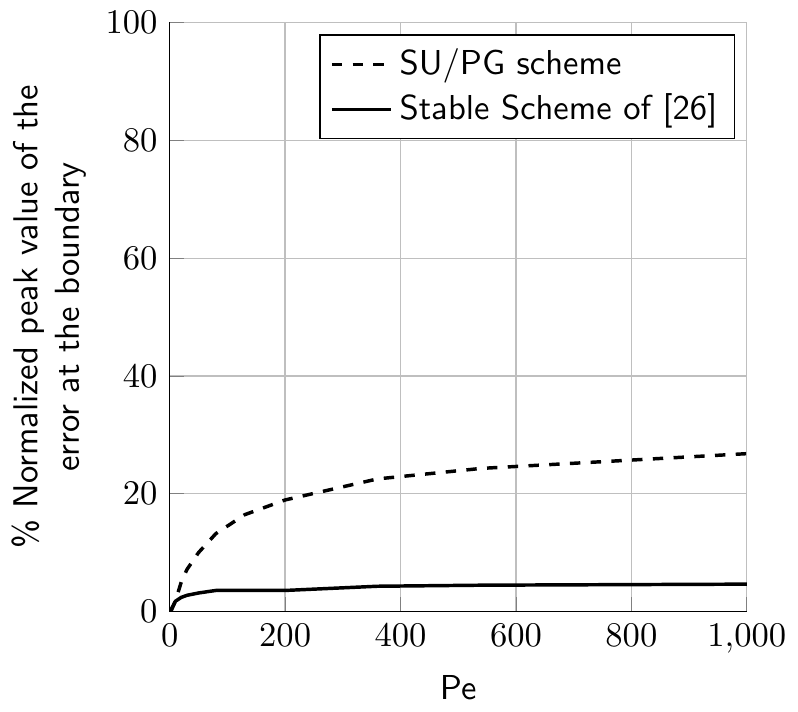}}}
		\caption{Selected results from the 2D problem. (a) $\nabla \cdot {\bf A} $ obtained from the SU/PG scheme (b) $\nabla \cdot {\bf A} $ obtained from Scheme presented in \cite{su1}  (c) Peak values $\nabla \cdot {\bf A}$ over a range of $Pe$ obtained from SU/PG scheme and the {stable scheme reported in} \cite{su1}.  (d) Peak values of the error over a range of $Pe$ {as} obtained from SU/PG scheme and the {stable scheme reported in} \cite{su1}.}
		\label{2dresults}
	\end{figure*}

During the analysis, it is found that the use of Coulomb's gauge i.e.,$\nabla \cdot {\bf A} = 0$ condition, is being grossly violated in the SU/PG scheme, which is not the case with the GFEM based schemes.
{As an example, a sample result of $\nabla \cdot \bf A$ computed from the SU/PG scheme is presented in the Fig. {\ref{divsA}}, where the presence of error could be clearly seen. Recently, a simple alternative scheme which is based on the classical GFEM has been presented in  }\cite{su1}   and it is shown to be very stable. 
The computed $\nabla \cdot \bf A$ obtained by employing this scheme is also presented, in Fig. {\ref{divaA}.}

\begin{figure}
		\centering
		\mbox{\subfloat[]{\label{supgbx} 
\includegraphics[scale=0.44]{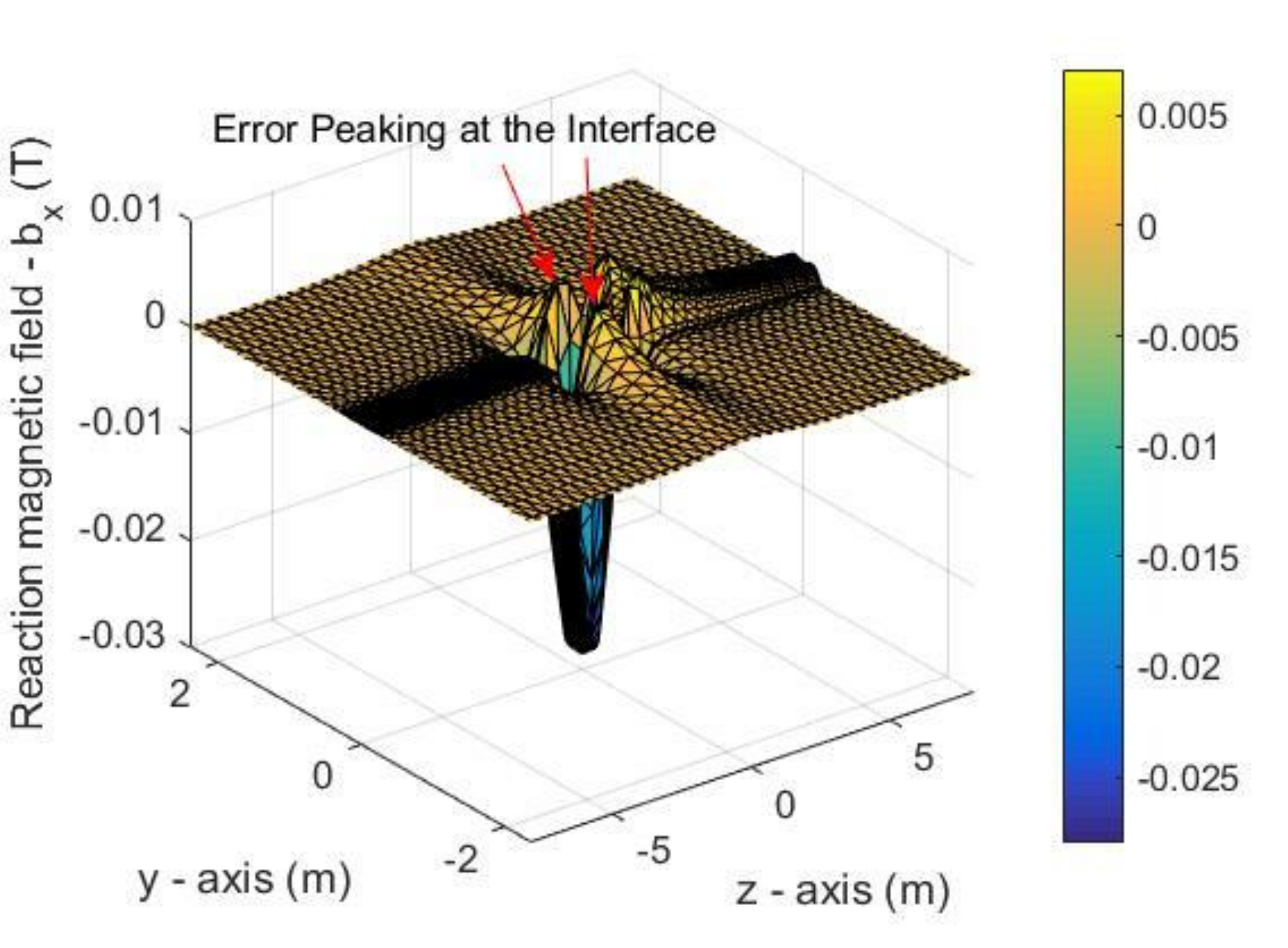}}}  \\
\mbox{\subfloat[]{\label{supgbxdivA} 
\includegraphics[scale=0.44]{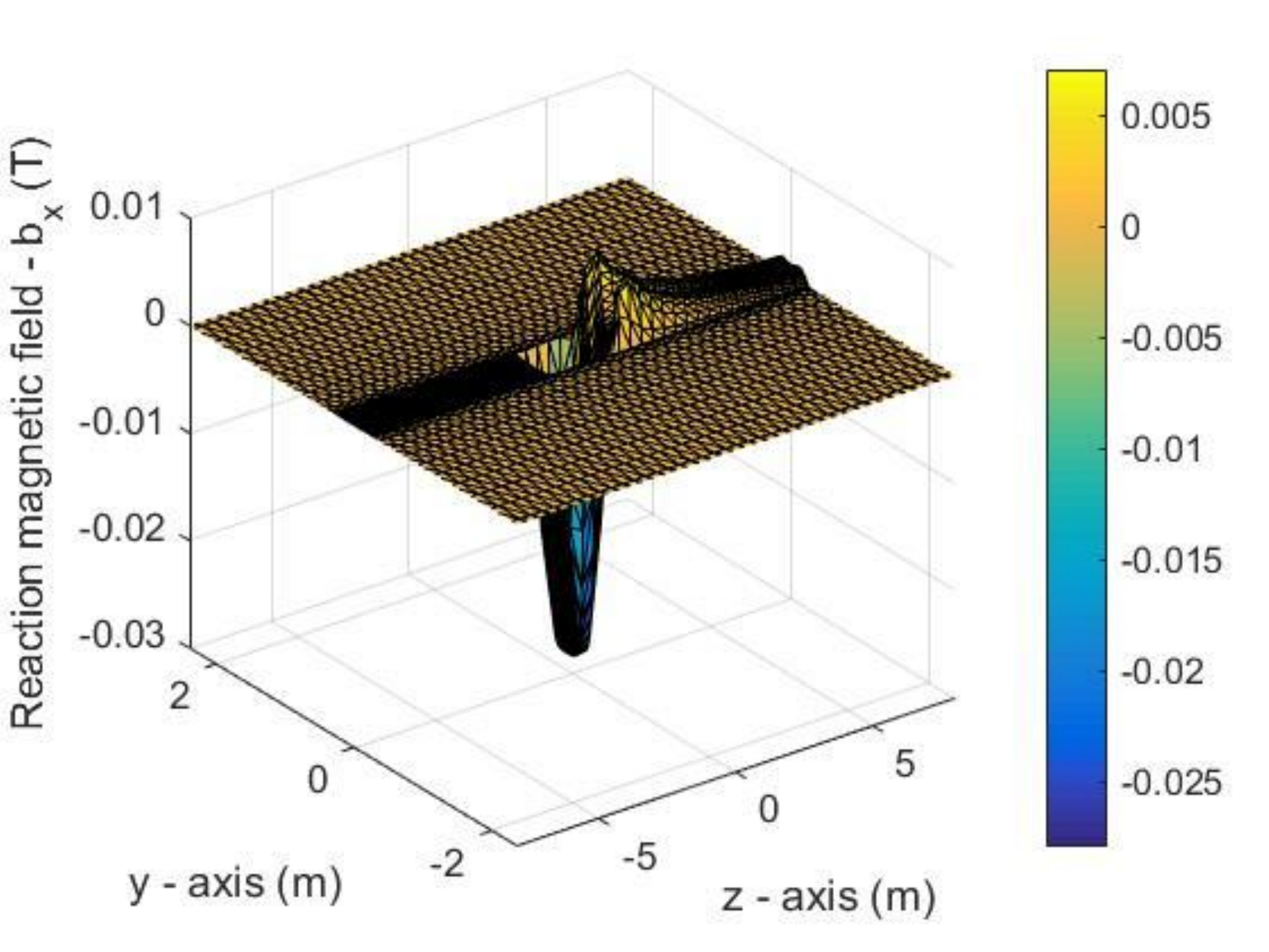}}} \\
\mbox{\subfloat[]{\label{bxref} 
\includegraphics[scale=0.44]{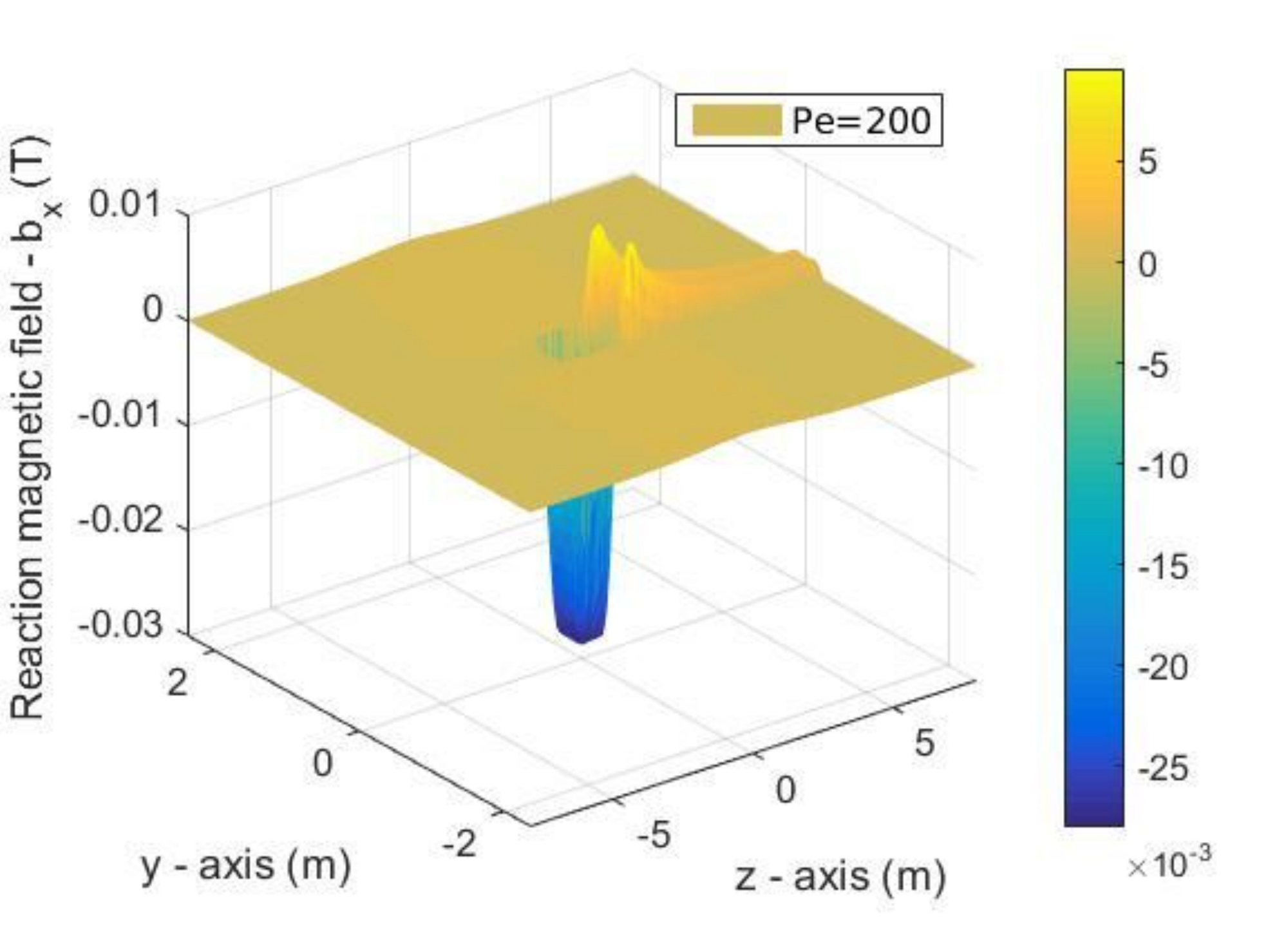}}}
		\caption{Reaction magnetic field $b_x$ obtained with the SU/PG scheme for $Pe=200$. (a) With the assumption of $\nabla \cdot {\bf A} = 0$  (b)  Without the assumption of $\nabla \cdot {\bf A} = 0$ (c) Reference solution generated with very fine discretisation \cite{nemosu}}
		\label{2dresults_supg}
	\end{figure}

In order to make a meaningful assessment, for each of the $Pe$, peak values of the $b_x$ is taken for normalisation. The required values are taken from the results presented in \cite{nemosu}, which has been reproduced as Fig. {\ref{bxref}} for a ready reference.
By comparing the results, it can be seen that the value of error in $\nabla \cdot \bf A$ is about 30\% in the SU/PG solution, while it is only about 5\% for the scheme proposed in \cite{su1}. This suggests a method for a direct correlation between the transverse boundary error and the peak magnitude of $\nabla \cdot {\bf A}$.
To get a more quantitative picture, both the peak error as well as the peak $\nabla \cdot {\bf A}$ are obtained for a range of $Pe$ and plotted in Figs. \ref{dAerr} and \ref{Bderr}.
 It is clear that the peak values of $\nabla \cdot {\bf A}$ and the peak boundary error with SU/PG scheme follows the exact same trend over the range of $Pe$.

This exercise clearly indicate that the boundary error in the solution is strongly correlated with the existence of $\nabla \cdot {\bf A} = 0$ in the weak solution through SU/PG scheme.
It is worth noting here that the Coulombs gauge is commonly employed in moving conductor literature \cite{mc5mc1, mc2av2, mc6tf1, mc5mc2, mc4ge1} especially for the nodal formulation. The following vector identity is employed most of the times,
\begin{equation}\label{vecId}
\nabla \times \nabla \times \bf A = \nabla (\nabla \cdot {\bf A}) - \nabla^2 \bf A
\end{equation}
assuming $\nabla \cdot \bf A = 0$,
\begin{equation}\label{divz}
\nabla \times \nabla \times \bf A =  - \nabla^2 \bf A
\end{equation}

However, the results of SU/PG scheme does not confirm to $\nabla \cdot {\bf A} = 0$ and therefore, it is necessary to avoid the use of Coulomb's gauge (refer to Figs. \ref{divsA} \& \ref{dAerr}). For this the corresponding term must be retained in the formulation.
Please note that this might result in non-unique value of $\bf A$, however, it is  $\nabla \times {\bf A}$ is the variable of interest not $\bf A$ directly.

{As the weak form for the Ampere's law (which involves $\nabla \times \nabla \times {\bf{A}}$) is not readily available for scalar elements, it will be dealt below.}

The weighted residue formulation for the $x$-component of  $\nabla \times \dfrac{1}{\mu} \nabla \times  \bf A$ is,
\begin{align}\label{aeq1}
 \int_\Omega N_g ~(\nabla \times \dfrac{1}{\mu} \nabla \times {\bf A} )\cdot\hat{\bf x} ~d\Omega 
\end{align}
where $N_g$ is the Galerkin weighting function and $\hat{\bf x}$ is the unit vector along the $x$-axis. Applying the vector identity given in (\ref{vecId}) on (\ref{aeq1}),
\begin{align} \label{aeq2}
\Rightarrow  \int_\Omega N_g ~\nabla \cdot\Big(  \dfrac{1}{\mu}  \dfrac{\partial {\bf A}}{\partial x}  - \dfrac{1}{\mu} \nabla  A_x \Big) ~d\Omega 
\end{align}
Applying integration by parts for (\ref{aeq2}),
\begin{equation} \label{aeq3}
\begin{split}
\Rightarrow~ - \int_\Omega  ~\nabla N_g \cdot\Big(  \dfrac{1}{\mu}  \dfrac{\partial {\bf A}}{\partial x}  - \dfrac{1}{\mu} \nabla  A_x \Big) ~d\Omega ~\dots \\[3mm]+ \oint_\Gamma N_g ~\Big(  \dfrac{1}{\mu}  \dfrac{\partial {\bf A}}{\partial x}  - \dfrac{1}{\mu} \nabla  A_x \Big) \cdot \hat{\bf n} ~d\Gamma 
\end{split}
\end{equation}
where, $\hat{\bf n}$ is the unit normal vector. Remember that in the Laplacian formulation of scalar potential, by neglecting the boundary integral term in the weak form, the inter element boundary condition ${\bf D}_1 \cdot {\hat{\bf n}} = {\bf D}_2 \cdot {\hat{\bf n}}$ or ${\bf J}_1 \cdot {\hat{\bf n}} = {\bf J}_2 \cdot {\hat{\bf n}}$ get automatically satisfied \cite{emfmbook}. Thus, it was not necessary to compute the boundary integral term there and that greatly improved the computational efficiency. Considering this, it is worthwhile to see the effect of neglecting the boundary integral term in the equation (\ref{aeq3}). The boundary integral term can be simplified as follows,

\begin{equation} \label{bieq1}
\begin{split}
\Rightarrow~ & \oint_\Gamma N_g ~\Big(  \dfrac{1}{\mu}  \dfrac{\partial {\bf A}}{\partial x}  - \dfrac{1}{\mu} \nabla  A_x \Big) \cdot \hat{\bf n} ~d\Gamma \\[2mm]
&=  \oint_\Gamma  \dfrac{N_g}{\mu} ~\Big(   \dfrac{\partial { A_y}}{\partial x} n_y + \dfrac{\partial { A_z}}{\partial x} n_z  -  \dfrac{\partial { A_x}}{\partial y} n_y - \dfrac{\partial { A_x}}{\partial z} n_z \Big)  ~d\Gamma  \\[2mm]
&=  \oint_\Gamma  {N_g} ~(h_z n_y - h_y n_z )  ~d\Gamma  
\end{split}
\end{equation}
where, \[ h_z = \dfrac{1}{\mu} \Big(  \dfrac{\partial { A_y}}{\partial x} -  \dfrac{\partial { A_x}}{\partial y}  \Big) ~~ \text{and}~~ h_y = \dfrac{1}{\mu} \Big(  \dfrac{\partial { A_x}}{\partial z} -  \dfrac{\partial { A_z}}{\partial x}  \Big) \]

Equation (\ref{bieq1}) corresponds to the  $x$-component, the boundary integral terms for other components can also be derived in the same fashion. The boundary integral term for the $y$-component is,

\begin{equation}\label{biyeq}
 \oint_\Gamma  {N_g} ~(h_x n_z - h_z n_x )  ~d\Gamma  
\end{equation}
and that for the $z$-component is,
\begin{equation}\label{bizeq}
 \oint_\Gamma  {N_g} ~(h_y n_x - h_x n_y )  ~d\Gamma  
\end{equation}
Combining (\ref{bieq1}), (\ref{biyeq}) and (\ref{bizeq}), the vector form of the boundary integral terms can be written as,
\begin{equation}\label{bieq_}
\Rightarrow ~- \oint_\Gamma  {N_g} ~({\bf h }\times {\bf n} )  ~d\Gamma  
\end{equation} 
Equating (\ref{bieq_}) to zero results in maintaining the tangential continuity of the magnetic field intensity $\bf h$ \cite{efem1}, i.e., 
\begin{equation}\label{hteq}
\bf h_{1} \cdot \hat{t} \approx h_{2} \cdot \hat{t}
\end{equation} 
across the elements 1 and 2. Thus, neglecting the boundary integral in (\ref{aeq3}), helps to maintain the continuity of magnetic field intensity across the element interface. Therefore, the weak formulation contains only the volume integral term in the proposed approach and the complete SU/PG-weak formulation of (\ref{eqge1}) and (\ref{eqge2}) are provided in the Appendix \ref{appA}.
\begin{figure}
		\centering
		\mbox{\subfloat[]{\label{supgj} 
\includegraphics[scale=0.44]{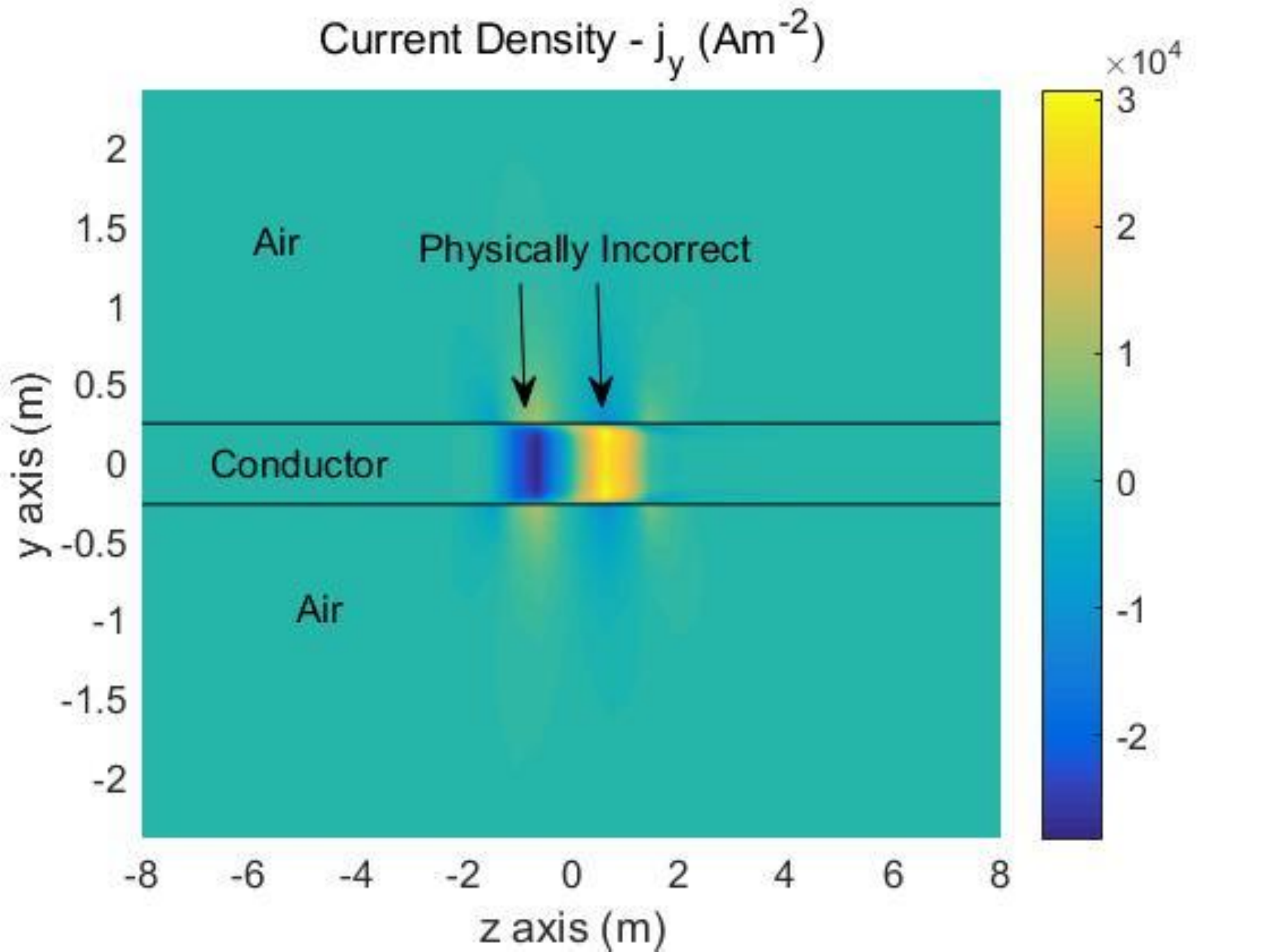}}}  \\
\mbox{\subfloat[]{\label{supgjdivA} 
\includegraphics[scale=0.44]{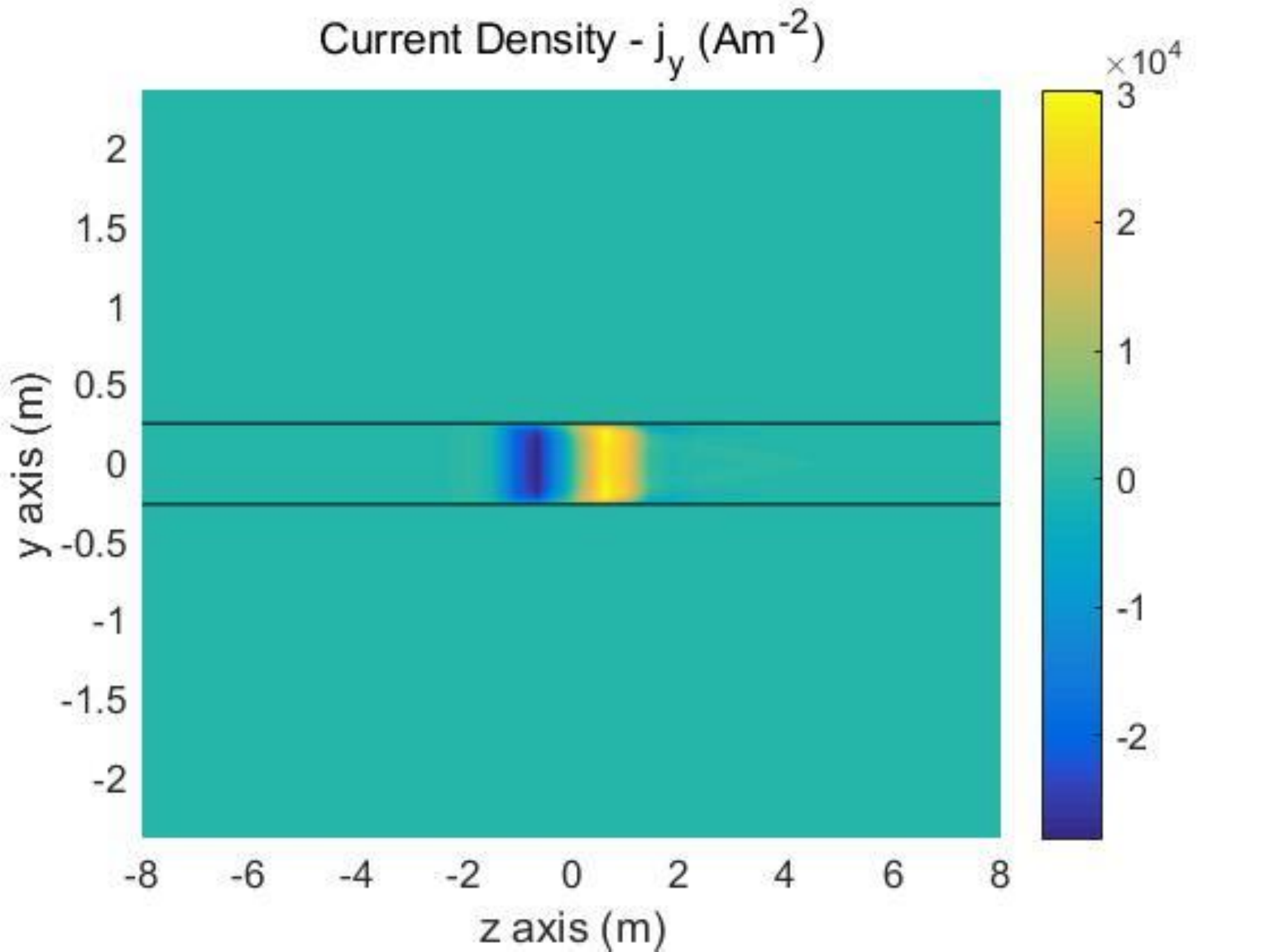}}}
		\caption{Current density $j_y = \partial h_x/\partial z$, obtained with the SU/PG scheme for $Pe=200$. (a) With the assumption of $\nabla \cdot {\bf A} = 0$  (b)  Without the assumption of $\nabla \cdot {\bf A} = 0$ }
		\label{2dresults_supg_j}
\end{figure}

\begin{figure}
		\centering
\includegraphics[scale=1]{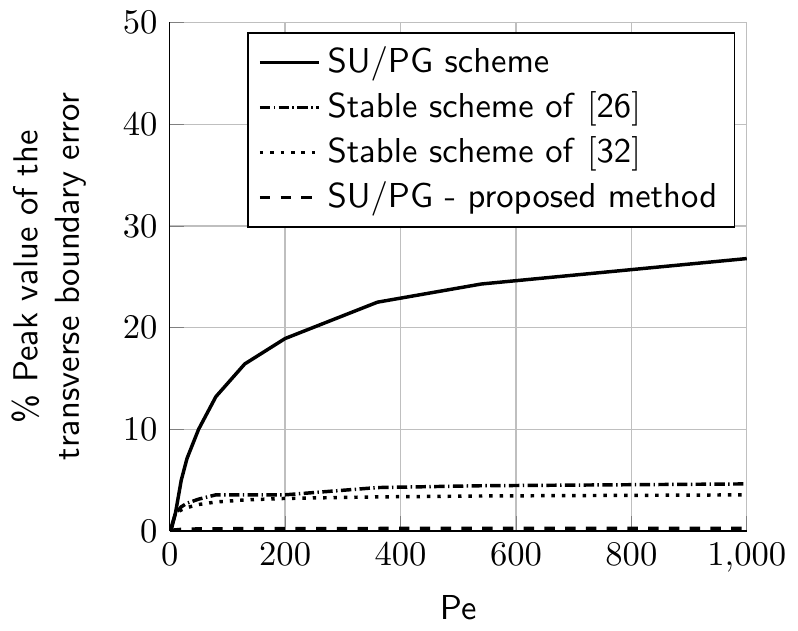}
		\caption{Comparison of \% peak value of transverse boundary error present in different schemes}
		\label{errprop} 
\end{figure}

With this, the simulations are carried out for  range of Peclet numbers {and the sample simulation results are presented in Fig. {\ref{2dresults_supg}}. When the $\nabla \cdot {\bf{A}}$ is retained with the SU/PG scheme, error at the transverse boundary virtually disappears.} Also, It can be clearly seen that the reaction magnetic field profile obtained with the proposed approach (Fig.  \ref{supgbxdivA}) is matching exactly with the reference solution presented in Fig. \ref{bxref}. Consequently, the non-physical solution of $\nabla \times {\bf h} \neq 0$ in the adjacent air-region, is also seem to be practically absent  with the proposed approach (refer to Fig. \ref{2dresults_supg_j}, {which provides solution for both}). { For a range of $Pe$ considered in Fig. {\ref{dAerr}}, the simulation results with $\nabla \cdot {\bf A}$ retained is found to have peak error $ < 1\%$ (refer to Fig {\ref{errprop}}).}
{
In the next subsection, the effect of avoiding Coulomb's gauge on the condition number of the FEM matrix is discussed.}

\subsection{Effect on condition number}

\newcolumntype{C}[1]{>{\centering\let\newline\\\arraybackslash\hspace{0pt}}m{#1}}
\begin{table*} [tbp]
\sffamily
\centering
\caption{Comparison of Condition number and accuracy}
\label{Tab:1}
{\renewcommand{\arraystretch}{1.6}
\begin{tabular}{@{}|C{0.1\textwidth}|C{0.1\textwidth}|C{0.1\textwidth}|C{0.1\textwidth}|C{0.1\textwidth}|C{0.1\textwidth}|C{0.1\textwidth}|C{0.1\textwidth}|}
\hline 
\multicolumn{1}{|C{0.1\textwidth}}{} \vline & \multicolumn{3}{C{0.3\textwidth}}{SU/PG scheme with Proposed method (Node element) }\vline & \multicolumn{2}{C{0.2\textwidth}}{Edge Element} \vline & \multicolumn{2}{C{0.2\textwidth}}{SUPG Scheme (Node element)} \vline \\ 
\cline{2-8} 
 Pe& $\alpha$ & Condition no. & \% Peak Error & Condition no. & \% Peak Error & Condition no. & \% Peak Error \\
\hline
 & 0.00 & 8.9415e+11 & 0.221 &  &  &  & \\
\cline{2-4}
50 & 0.05 & 4.1364e+08 & 2.215 & 4.8747e+09 & 22.516 & 1.7352e+08& 10.007\\
\cline{2-4}
 & 0.10 & 3.2260e+08 & 3.466 &  & & & \\
 \hline   
 & 0.00 & 1.5178e+12 & 0.257 &  &  &  & \\
\cline{2-4} 
200  & 0.05 & 2.0991e+09 & 3.359 & 1.9438e+10 & 25.018&1.8165e+09 & 18.942\\
\cline{2-4}
  & 0.10 & 1.9721e+09 & 5.361 &  & & & \\
 \hline
 & 0.00 & 2.1925e+14 & 0.271 &  &  &  & \\
\cline{2-4}
1000 & 0.05 & 4.0514e+11 & 4.289 &3.3003e+11  &25.733&1.8677e+11& 30.021\\
\cline{2-4}
 & 0.10 & 3.4614e+11 & 6.791 &  & & & \\
 \hline
\end{tabular}}
\end{table*} 

It is evident that the avoiding Coulomb's gauge from the governing equation can impact the condition number of the final FEM matrix. In order to address this, a possible remedy is discussed in this subsection. In this, a scaled amount of $\nabla \cdot \bf{A}$ is injected into the original governing equation, with the scaling parameter $\alpha$. The results of this new variation is tabulated in Table {\ref{Tab:1}}. 
For comparison, simulations also carried using the edge elements and the results of same is included in Table {\ref{Tab:1}}. As mentioned in the introduction, the edge elements cannot inherently avoid the erroneous results arising out of the convection dominated simulations.

As shown in the table {\ref{Tab:1}}, the \% peak error in the proposed scheme can be seen to slowly increase as we increase the $\nabla \cdot \bf{A}$ contribution. However, for the values of $\alpha$ listed in the table, the accuracy of the proposed scheme is substantially better than the edge element method or the SU/PG scheme.  When $\alpha = 0.05$ or in other words, when the $\nabla \cdot \bf{A}$ addition is around 5\%, the accuracy at the material interface is at least an order of magnitude better than the other methods; while keeping the condition number in a comparably reasonable value.

It can also be observed that the condition number from the edge element method is not better than the nodal element method. This may be attributed to the absence of normal continuity across the edge elements \cite{r:edgeadv, r:edgefall}.

\subsection{Verification with 3D electromagnetic flowmeter problem}
\begin{figure}
		\centering
\includegraphics[scale=0.44]{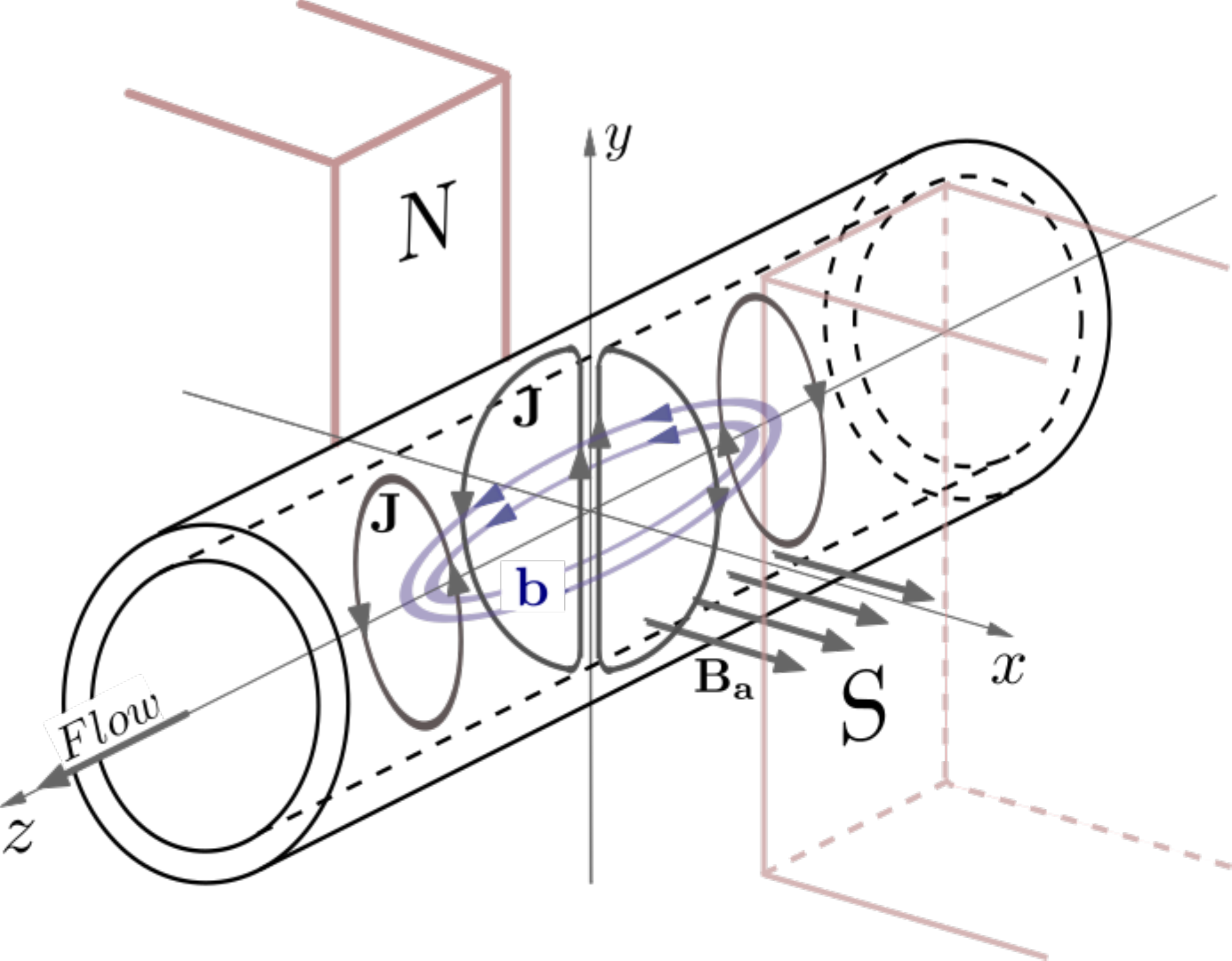}
		\caption{Schematic of the electromagnetic flowmeter, shown with the
		current density ${\bf J}$ and the reaction magnetic field ${\bf b}$}
		\label{f:emfm}
\end{figure}

\begin{figure}
		\centering
		\mbox{\subfloat[]{\label{supgj3d} 
\includegraphics[scale=0.44]{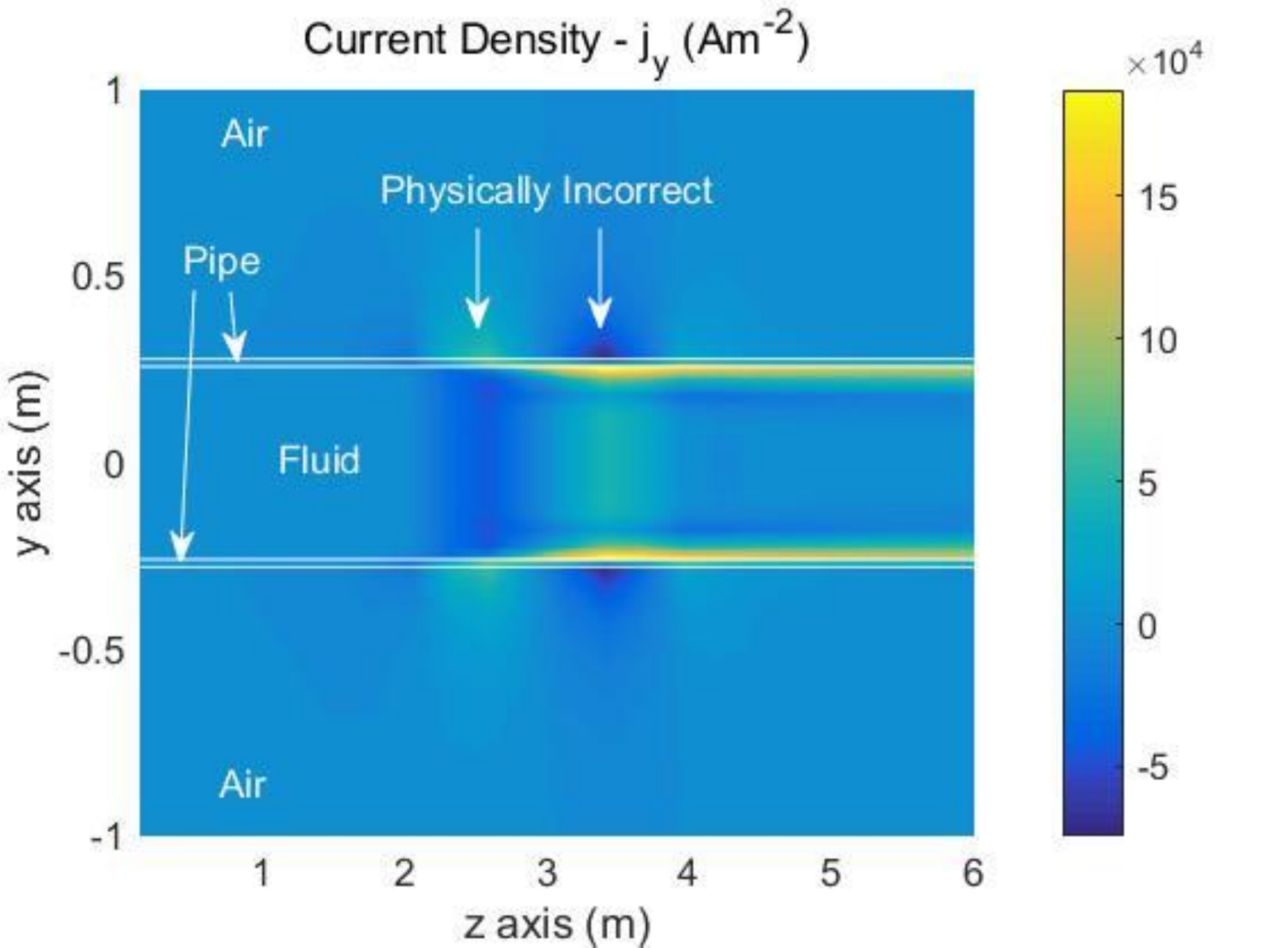}}}  \\
\mbox{\subfloat[]{\label{supgjdivA3d} 
\includegraphics[scale=0.44]{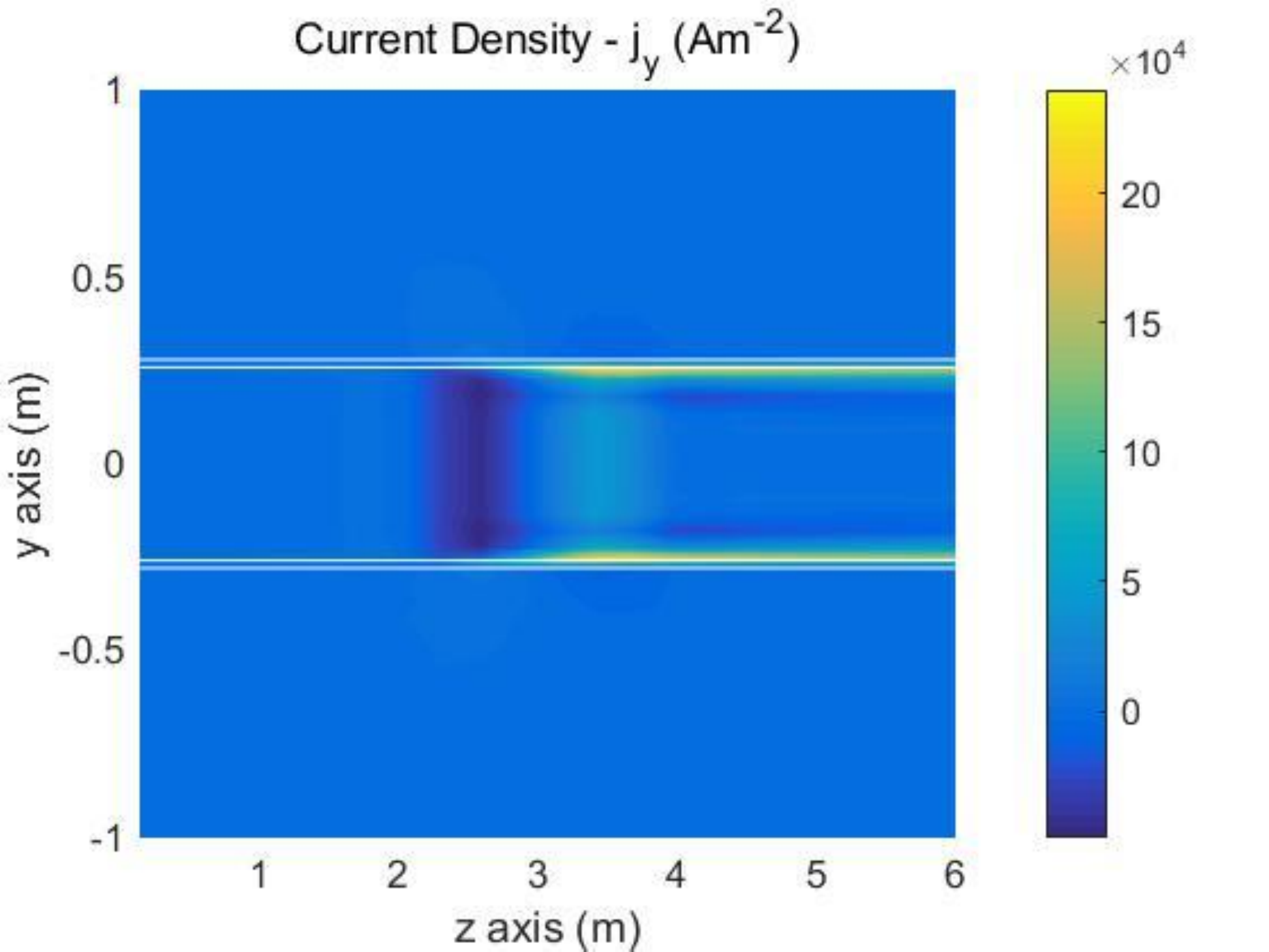}}}
		\caption{Current density $j_y = \partial h_x/\partial z - \partial h_z/\partial x $ in $x=0$ plane, obtained with the SU/PG scheme for $Pe=1600$. (a) With the assumption of $\nabla \cdot {\bf A} = 0$  (b)  Without the assumption of $\nabla \cdot {\bf A} = 0$ }
		\label{3dresults_supg_j}
\end{figure}

As a demonstration of the efficacy of the proposed approach in real life situations, a practical 3D problem, the electromagnetic flowmeter, is simulated for high Peclet numbers. The simulation is carried out using the bilinear nodal elements. The simulation parameters and the input quantities are taken as mentioned in \cite{su1} and the schematic diagram of the electromagnetic flowmeter is shown in Fig. \ref{f:emfm}, with its $x, y, z-axes$ defined. The simulations are carried out for a range of Peclet numbers and the results are confirming the conclusions, which are drawn for the 2D problem. The sample simulation results are presented in Fig. \ref{3dresults_supg_j}, with $Pe=1600$. From Fig. \ref{supgj3d} the non-physical current distribution in the air region can be seen, which is obtained from the SU/PG scheme, on the other hand, {the retention of $\nabla \cdot \bf{A}$ in the SU/PG scheme} provides a physically confirming solution (refer to Fig. \ref{supgjdivA3d}).

\begin{figure}
		\centering
\includegraphics[scale=0.44]{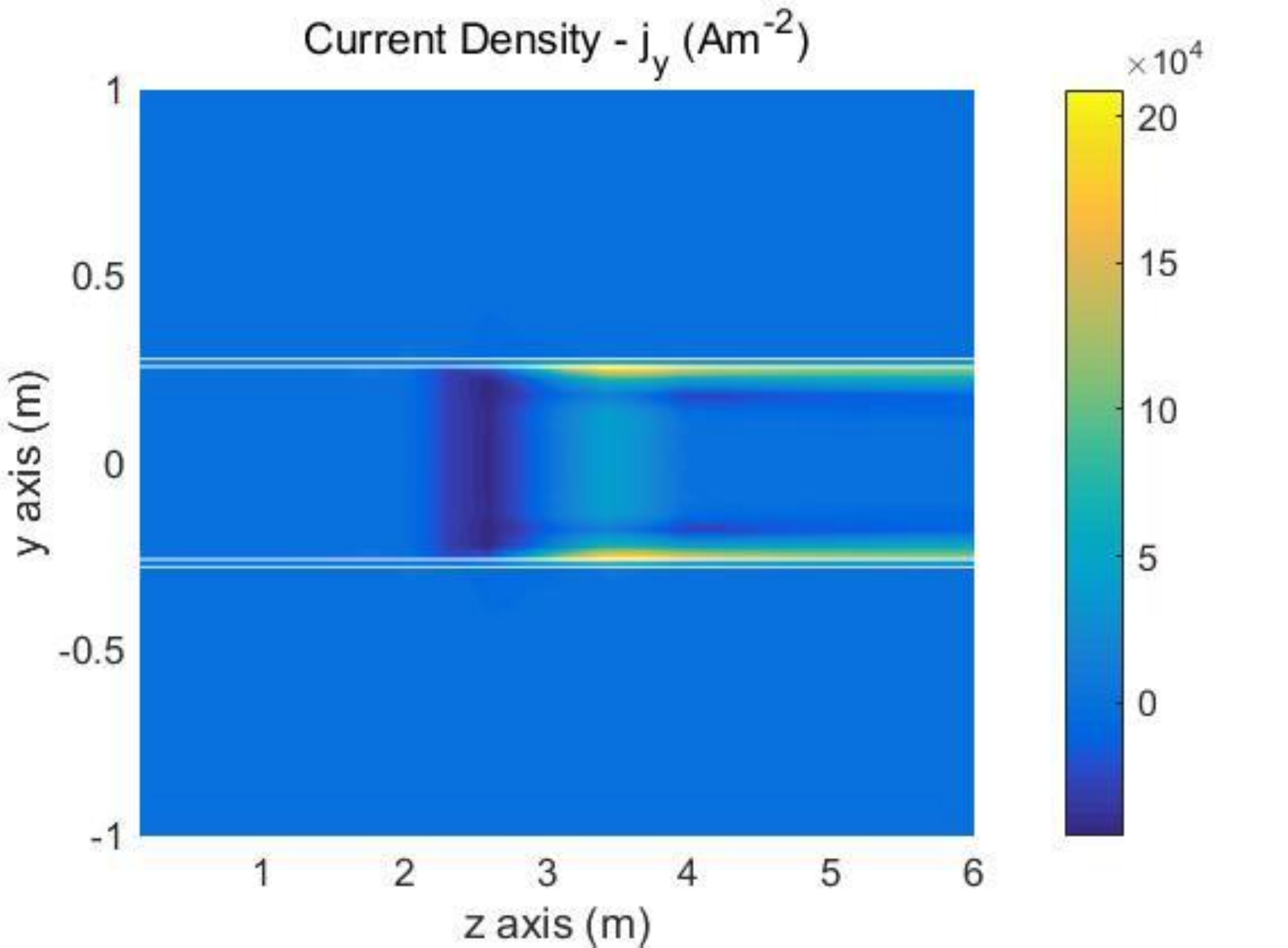}
		\caption{Current density $j_y = \partial h_x/\partial z - \partial h_z/\partial x $ in $x=0$ plane, obtained with the pole-zero cancellation based scheme \cite{su2} for $Pe=1600$.}
		\label{3dresults_schb_j}
\end{figure}
Incidentally, the recent pole-zero cancellation based, stabilization schemes proposed in \cite{su1, su2} observed to give accurate solutions without the boundary error. A sample simulation result for $Pe=1600$ with the scheme proposed in \cite{su2}, is provided in Fig. \ref{3dresults_schb_j}.

Referring to Fig. \ref{f:emfm}, along the $zy$-plane the complete circulation of induced current is represented. In the orthogonal direction, i.e., along
the $zx$-plane the currents are orthogonal to the plane and they are mostly
represented by the $y$-component of the current density $j_y$. The currents
along the $zx$-plane does not get directly oriented by the boundary, unlike the circulation along the $zy$-plane. Owing to this, the boundary error is expected
to be minimal along the $zx$-plane and the same is observed with this 3D
simulation. To confirm this observation further, the 2D version of the
moving conductor problem along the $zx$-plane is simulated. The results from this
simulation indicate insignificant boundary error, along the $zx$-plane. 

\section{Summary and Conclusion}
The SU/PG scheme is widely adopted upwinding technique, for high flow velocity problems due to its simplicity and being void of crosswind of diffusion \cite{quada1}. Besides, the transverse boundary error present in the SU/PG scheme is a serious issue and several non-linear solutions have been proposed (called as SOLD schemes) by the computational fluid dynamics community. Nevertheless, they suffer from computational inefficiency, as well as, some inaccuracy  \cite{soldreview1, soldreview2}. In addition, these techniques are rarely employed for moving conductor problems.

In the present work, a 2D version of the moving conductor problem, is taken for the analysis. The 2D problem represents the complete circulation of current with in the conducting region and thus it consists of 2-components of vector potential. Here, the source of the transverse boundary error is identified {with the enforcement of Coulomb gauge $\nabla \cdot{\bf A} = 0$. Therefore, the $\nabla \cdot {\bf{A}}$ term is retained and the weak form of the Amperes' law has been derived for scalar elements. Then, through the numerical simulation of 2D, as well as, 3D problems, it is shown that the retention of $\nabla \cdot {\bf{A}}$ practically eliminates the material interface boundary error present in the SU/PG scheme. } In addition, this approach preserves the linearity, as well as, computational efficiency of the SU/PG scheme.


%

\appendices
\section{Weak Form of the Governing Equations}\label{appA}
The weak form of (\ref{eqge2}) can be commonly found in the computational electromagnetic literature. Nevertheless, for the sake of completeness and quick reference, it has been provided below, along with the weak form of (\ref{eqge1}), which is derived.
The weak form of (\ref{eqge2}) can be written as,
\begin{equation}\label{w_eqge2}
\begin{split}
\int_\Omega  \sigma \nabla N_g \cdot  \nabla \phi ~d\Omega - \int_ \Omega \sigma \nabla N_g  \cdot ({\bf u} \times \nabla \times {\bf{A}})~d\Omega ~\dots \\ =  \int_ \Omega \sigma \nabla N_g  \cdot   ({\bf u} \times {\bf B_a})~d\Omega
\end{split}
\end{equation}
The advection-diffusion form of (\ref{eqge1}) requires the modification of the Galerkin weighting function $N_g$ to the SU/PG weighting function $N_s$, which is given by, 
\begin{equation}\label{supg_ns}
N_s  = N_g + \tau \dfrac{ \Delta l~{\bf u} }{ 2 \| {\bf u} \|} \cdot \nabla N_g
\end{equation}
where $\tau$ is the stabilization parameter, defined as $\tau \equiv coth(|Pe|) - {1}/{|Pe|}$. The stabilization parameter $\tau$ decides the optimal (right) amount of upwinding which is being added to the Galerkin weight function $N_g$. This added term, is employed to only the element interiors and in case of linear elements, the second derivative vanishes and the modified weighting function $N_s$ is applied only to the advective/first derivative terms \cite{multiscale1, supg1, quada1}. Therefore, the weak form of the $x$-component of (\ref{eqge1}) takes the form:
\begin{equation}\label{w_eqge1x}
\begin{split}
\int_\Omega \sigma N_s \dfrac{\partial \phi}{\partial x} ~d\Omega~- \int_\Omega  ~\nabla N_g \cdot\Big(  \dfrac{1}{\mu}  \dfrac{\partial {\bf A}}{\partial x}  - \dfrac{1}{\mu} \nabla  A_x \Big) ~d\Omega ~\dots \\[3mm]- \int_\Omega \sigma N_s~\Big( u_y  \dfrac{\partial A_y}{\partial x} - u_y  \dfrac{\partial A_x}{\partial y} - u_z  \dfrac{\partial A_x}{\partial z} + u_z  \dfrac{\partial A_z}{\partial x} \Big)~d\Omega ~\dots \\[3mm] 
+~ \alpha \int_\Omega ~\dfrac{\partial N_g}{\partial x} \nabla {\bf \cdot A} ~d\Omega
= \int_\Omega \sigma N_s (u_y B_z - u_z B_y)~ d\Omega
\end{split}
\end{equation}
and its $y$-component:
\begin{equation}\label{w_eqge1y}
\begin{split}
\int_\Omega \sigma N_s \dfrac{\partial \phi}{\partial y} ~d\Omega~- \int_\Omega  ~\nabla N_g \cdot\Big(  \dfrac{1}{\mu}  \dfrac{\partial {\bf A}}{\partial y}  - \dfrac{1}{\mu} \nabla  A_y \Big) ~d\Omega ~\dots \\[3mm]- \int_\Omega \sigma N_s~\Big( u_z  \dfrac{\partial A_z}{\partial y} - u_z  \dfrac{\partial A_y}{\partial z} - u_x  \dfrac{\partial A_y}{\partial x} + u_x  \dfrac{\partial A_x}{\partial y} \Big)~d\Omega ~\dots \\[3mm] 
+~ \alpha \int_\Omega ~\dfrac{\partial N_g}{\partial y} \nabla {\bf \cdot A} ~d\Omega
= \int_\Omega \sigma N_s (u_z B_x - u_x B_z) ~d\Omega
\end{split}
\end{equation}
and its $z$-component:
\begin{equation}\label{w_eqge1z}
\begin{split}
\int_\Omega \sigma N_s \dfrac{\partial \phi}{\partial z} ~d\Omega~- \int_\Omega  ~\nabla N_g \cdot\Big(  \dfrac{1}{\mu}  \dfrac{\partial {\bf A}}{\partial z}  - \dfrac{1}{\mu} \nabla  A_z \Big) ~d\Omega ~\dots \\[3mm]- \int_\Omega \sigma N_s~\Big( u_x  \dfrac{\partial A_x}{\partial z} - u_x  \dfrac{\partial A_z}{\partial x} - u_y  \dfrac{\partial A_z}{\partial y} + u_y  \dfrac{\partial A_y}{\partial z} \Big)~d\Omega ~\dots \\[3mm]
+~ \alpha \int_\Omega ~\dfrac{\partial N_g}{\partial z} \nabla {\bf \cdot A}~d\Omega
= \int_\Omega \sigma N_s (u_x B_y - u_y B_x ) ~d\Omega
\end{split}
\end{equation}

For the 2D problem considered, the $x$-component of the equation vanishes and only the $y$ and $z$ components remain. Also, the velocity is assumed to be along the $z$ direction (${\bf u} = u_z {\bf \hat{z}}$) and the input magnetic field is  $x$ directed (${\bf B_a} = B_x {\bf \hat{x}}$). With this, the weak forms of the  governing equations (\ref{eqge1}) and (\ref{eqge2}) are respectively given by,
\begin{equation}\label{w_eqge1y_2d}
\begin{split}
\int_\Omega \sigma N_s \dfrac{\partial \phi}{\partial y} ~d\Omega~+ \int_\Omega  ~ \dfrac{1}{\mu}  \Big( \dfrac{\partial  N_g}{\partial z} \dfrac{\partial  A_y}{\partial z} -\dfrac{\partial  N_g}{\partial z} \dfrac{\partial  A_z}{\partial y}   \Big) ~d\Omega ~\dots \\[3mm]
+ \int_\Omega \sigma N_s \Big( u_z  \dfrac{\partial A_y}{\partial z} - u_z  \dfrac{\partial A_z}{\partial y} \Big)~d\Omega ~\dots \\[3mm]
+~ \alpha \int_\Omega ~\dfrac{\partial N_g}{\partial y} 
\Big( \dfrac{\partial A_y}{\partial y} + \dfrac{\partial A_z}{\partial z} \Big)
~d\Omega 
 = \int_\Omega \sigma N_s u_z B_x ~d\Omega
\end{split}
\end{equation}
\begin{equation}\label{w_eqge1z_2d}
\begin{split}
\int_\Omega \sigma N_s \dfrac{\partial \phi}{\partial z} ~d\Omega~+ \int_\Omega  ~ \dfrac{1}{\mu}  \Big( \dfrac{\partial  N_g}{\partial y} \dfrac{\partial  A_z}{\partial y} -\dfrac{\partial  N_g}{\partial y} \dfrac{\partial  A_y}{\partial z}   \Big) ~d\Omega  ~\dots \\[3mm]
+~ \alpha \int_\Omega ~\dfrac{\partial N_g}{\partial z} 
\Big( \dfrac{\partial A_y}{\partial y} + \dfrac{\partial A_z}{\partial z} \Big)
~d\Omega 
 = 0
\end{split}
\end{equation}
\begin{equation}\label{w_eqge2}
\begin{split}
\int_\Omega  \sigma   \Big( \dfrac{\partial  N_g}{\partial y} \dfrac{\partial \phi}{\partial y} +\dfrac{\partial  N_g}{\partial z} \dfrac{\partial \phi}{\partial z}   \Big)  ~d\Omega + \int_ \Omega \sigma u_z  \dfrac{\partial  N_g}{\partial y} \dfrac{\partial A_y}{\partial z}~ d\Omega ~\dots  \\[3mm] - \int_ \Omega \sigma u_z \dfrac{\partial  N_g}{\partial y} \dfrac{\partial A_z}{\partial y}  ~d\Omega  =  \int_ \Omega \sigma u_z B_x \dfrac{\partial  N_g}{\partial y}  ~d\Omega
\end{split}
\end{equation}




\ifCLASSOPTIONcaptionsoff
  \newpage
\fi



%
\bibliographystyle{IEEEtran} 
\bibliography{./References/References}

\begin{thebibliography}{10}
\providecommand{\url}[1]{#1}
\csname url@samestyle\endcsname
\providecommand{\newblock}{\relax}
\providecommand{\bibinfo}[2]{#2}
\providecommand{\BIBentrySTDinterwordspacing}{\spaceskip=0pt\relax}
\providecommand{\BIBentryALTinterwordstretchfactor}{4}
\providecommand{\BIBentryALTinterwordspacing}{\spaceskip=\fontdimen2\font plus
\BIBentryALTinterwordstretchfactor\fontdimen3\font minus
  \fontdimen4\font\relax}
\providecommand{\BIBforeignlanguage}[2]{{%
\expandafter\ifx\csname l@#1\endcsname\relax
\typeout{** WARNING: IEEEtran.bst: No hyphenation pattern has been}%
\typeout{** loaded for the language `#1'. Using the pattern for}%
\typeout{** the default language instead.}%
\else
\language=\csname l@#1\endcsname
\fi
#2}}
\providecommand{\BIBdecl}{\relax}
\BIBdecl

\bibitem{reducedb}
O.~Biro, K.~Preis, W.~Renhart, K.~Richter, and G.~Vrisk, ``Performance of
  different vector potential formulations in solving multiply connected 3-d
  eddy current problems,'' \emph{IEEE Trans. Magnetics}, vol.~26, no.~2, pp.
  438--441, 1990.

\bibitem{fmbase}
T.~Shimizu, N.~Takeshima, and N.~Jimbo, ``A numerical study on faraday-type
  electromagnetic flowmeter in liquid metal system, (i),'' \emph{Journal of
  Nuclear Science and Technology}, vol.~37, no.~12, pp. 1038--1048, 2000.

\bibitem{cdbook}
O.~Zienkiewicz, R.~Taylor, and P.~Nithiarasu, \emph{The Finite Element Method
  for Fluid Dynamics}.\hskip 1em plus 0.5em minus 0.4em\relax Elsevier Science,
  2005.

\bibitem{upfdm1_sp}
D.~Spalding, ``A novel finite difference formulation for differential
  expressions involving both first and second derivatives,''
  \emph{International Journal for Numerical Methods in Engineering}, vol.~4,
  no.~4, pp. 551--559, 1972.

\bibitem{upfdm1_ru}
A.~Runchal, ``Convergence and accuracy of three finite difference schemes for a
  two-dimensional conduction and convection problem,'' \emph{International
  Journal for Numerical Methods in Engineering}, vol.~4, no.~4, pp. 541--550,
  1972.

\bibitem{up1}
I.~Christie, D.~F. Griffiths, A.~R. Mitchell, and O.~C. Zienkiewicz, ``Finite
  element methods for second order differential equations with significant
  first derivatives,'' \emph{International Journal for Numerical Methods in
  Engineering}, vol.~10, no.~6, pp. 1389--1396, 1976.

\bibitem{mc6tf1}
M.~Odamura, ``Upwind finite element solution for saturated traveling magnetic
  field problems,'' \emph{Electrical Engineering in Japan}, vol. 105, no.~4,
  pp. 126--132, 1985.

\bibitem{mc6tf2}
M.~Ito, T.~Takahashi, and M.~Odamura, ``Up-wind finite element solution of
  travelling magnetic field problems,'' \emph{IEEE Trans. Magnetics}, vol.~28,
  no.~2, pp. 1605--1610, 1992.

\bibitem{mc2av1}
D.~Rodger, P.~Leonard, and T.~Karaguler, ``An optimal formulation for 3d moving
  conductor eddy current problems with smooth rotors,'' \emph{IEEE Trans.
  Magnetics}, vol.~26, no.~5, pp. 2359--2363, 1990.

\bibitem{mc3eb1}
E.~Chan and S.~Williamson, ``Factors influencing the need for upwinding in
  two-dimensional field calculation,'' \emph{IEEE Trans. Magnetics}, vol.~28,
  no.~2, pp. 1611--1614, 1992.

\bibitem{mc4ge1}
N.~Allen, D.~Rodger, P.~Coles, S.~Strret, and P.~Leonard, ``Towards increased
  speed computations in 3d moving eddy current finite element modelling,''
  \emph{IEEE Trans. Magnetics}, vol.~31, no.~6, pp. 3524--3526, 1995.

\bibitem{mc5mc1}
D.~Rodger, T.~Karguler, and P.~Leonard, ``A formulation for 3d moving conductor
  eddy current problems,'' \emph{IEEE Trans. Magnetics}, vol.~25, no.~5, pp.
  4147--4149, 1989.

\bibitem{mc2av2}
J.~Bird, T.~Lipo \emph{et~al.}, ``A 3-d magnetic charge finite-element model of
  an electrodynamic wheel,'' \emph{IEEE Trans. Magnetics}, vol.~44, no.~2, pp.
  253--265, 2008.

\bibitem{mc3eb2}
H.~Vande~Sande, H.~De~Gersem, and K.~Hameyer, ``Finite element stabilization
  techniques for convection-diffusion problems,'' \emph{7th International
  journal of theoretical electrotechnics}, pp. 56--59, 1999.

\bibitem{mcsupg_mcfit}
L.~Codecasa and P.~Alotto, ``2-d stabilized fit formulation for eddy-current
  problems in moving conductors,'' \emph{IEEE Trans. Magnetics}, vol.~51,
  no.~3, pp. 1--4, 2015.

\bibitem{mcsupg_cable}
Y.~Liang, ``Steady-state thermal analysis of power cable systems in ducts using
  streamline-upwind/petrov-galerkin finite element method,'' \emph{Dielectrics
  and Electrical Insulation, IEEE Transactions on}, vol.~19, no.~1, pp.
  283--290, 2012.

\bibitem{mcsupg_mfluid}
S.~Noguchi and S.~Kim, ``Magnetic field and fluid flow computation of plural
  kinds of magnetic particles for magnetic separation,'' \emph{IEEE Trans.
  Magnetics}, vol.~48, no.~2, pp. 523--526, 2012.

\bibitem{edgeup2}
E.~X. Xu, J.~Simkin, and S.~C. Taylor, ``Streamline upwinding in a 3-d
  edge-element method modeling eddy currents in moving conductors,'' \emph{IEEE
  transactions on magnetics}, vol.~42, no.~4, pp. 667--670, 2006.

\bibitem{r:edgeup3}
F.~Henrotte, H.~Heumann, E.~Lange, and K.~Hameyer, ``Upwind 3-d vector
  potential formulation for electromagnetic braking simulations,'' \emph{IEEE
  Transactions on Magnetics}, vol.~46, no.~8, pp. 2835--2838, 2010.

\bibitem{discop}
T.~J. Hughes, M.~Mallet, and M.~Akira, ``A new finite element formulation for
  computational fluid dynamics: Ii. beyond supg,'' \emph{Computer Methods in
  Applied Mechanics and Engineering}, vol.~54, no.~3, pp. 341--355, 1986.

\bibitem{nemosu}
S.~Subramanian and U.~Kumar, ``Existence of boundary error transverse to the
  velocity in su/pg solution of moving conductor problem,'' in \emph{Numerical
  Electromagnetic and Multiphysics Modeling and Optimization (NEMO), 2016 IEEE
  MTT-S International Conference on}.\hskip 1em plus 0.5em minus 0.4em\relax
  IEEE, 2016, pp. 1--2.

\bibitem{discop1}
R.~Codina, ``A discontinuity-capturing crosswind-dissipation for the finite
  element solution of the convection-diffusion equation,'' \emph{Computer
  Methods in Applied Mechanics and Engineering}, vol. 110, no.~3, pp. 325--342,
  1993.

\bibitem{fic2}
E.~O{\~n}ate, F.~Z{\'a}rate, and S.~R. Idelsohn, ``Finite element formulation
  for convective--diffusive problems with sharp gradients using finite
  calculus,'' \emph{Computer methods in applied mechanics and engineering},
  vol. 195, no.~13, pp. 1793--1825, 2006.

\bibitem{soldreview1}
V.~John and P.~Knobloch, ``On spurious oscillations at layers diminishing
  (sold) methods for convection--diffusion equations: Part i--a review,''
  \emph{Computer Methods in Applied Mechanics and Engineering}, vol. 196,
  no.~17, pp. 2197--2215, 2007.

\bibitem{soldreview2}
------, ``On spurious oscillations at layers diminishing (sold) methods for
  convection--diffusion equations: Part ii--analysis for p1 and q1 finite
  elements,'' \emph{Computer Methods in Applied Mechanics and Engineering},
  vol. 197, no.~21, pp. 1997--2014, 2008.

\bibitem{su1}
S.~Subramanian and U.~Kumar, ``Augmenting numerical stability of the galerkin
  finite element formulation for electromagnetic flowmeter analysis,''
  \emph{IET Science, Measurement \& Technology}, vol.~10, no.~4, pp. 288--295,
  2016.

\bibitem{mc5mc2}
S.-K. Hong and N.~Ida, ``Modeling of velocity terms in 3d eddy current
  problems,'' \emph{IEEE Trans. Magnetics}, vol.~28, no.~2, pp. 1178--1181,
  1992.

\bibitem{emfmbook}
J.~Shercliff, \emph{The Theory of Electromagnetic Flow-Measurement}, ser.
  Cambridge Science Classics.\hskip 1em plus 0.5em minus 0.4em\relax Cambridge
  University Press, 1987.

\bibitem{efem1}
P.~P. Silvester and R.~L. Ferrari, \emph{Finite elements for electrical
  engineers}.\hskip 1em plus 0.5em minus 0.4em\relax Cambridge university
  press, 1996.

\bibitem{r:edgeadv}
G.~Mur, ``Edge elements, their advantages and their disadvantages,'' \emph{IEEE
  transactions on magnetics}, vol.~30, no.~5, pp. 3552--3557, 1994.

\bibitem{r:edgefall}
------, ``The fallacy of edge elements,'' \emph{IEEE Transactions on
  Magnetics}, vol.~34, no.~5, pp. 3244--3247, 1998.

\bibitem{su2}
S.~Subramanian and U.~Kumar, ``Stable galerkin finite-element scheme for the
  simulation of problems involving conductors moving rectilinearly in magnetic
  fields,'' \emph{IET Science, Measurement \& Technology}, vol.~10, no.~8, pp.
  952--962, 2016.

\bibitem{quada1}
T.-P. Fries and H.~G. Matthies, ``A review of petrov--galerkin stabilization
  approaches and an extension to meshfree methods,'' \emph{Technische
  Universitat Braunschweig, Brunswick}, 2004.

\bibitem{multiscale1}
T.~J. Hughes, ``Multiscale phenomena: Green's functions, the
  dirichlet-to-neumann formulation, subgrid scale models, bubbles and the
  origins of stabilized methods,'' \emph{Computer methods in applied mechanics
  and engineering}, vol. 127, no.~1, pp. 387--401, 1995.

\bibitem{supg1}
A.~N. Brooks and T.~J. Hughes, ``Streamline upwind/petrov-galerkin formulations
  for convection dominated flows with particular emphasis on the incompressible
  navier-stokes equations,'' \emph{Computer methods in applied mechanics and
  engineering}, vol.~32, no.~1, pp. 199--259, 1982.

\end{thebibliography}
%

%
\begin{IEEEbiography}[{\includegraphics[width=1in,height=1.25in,clip,keepaspectratio]{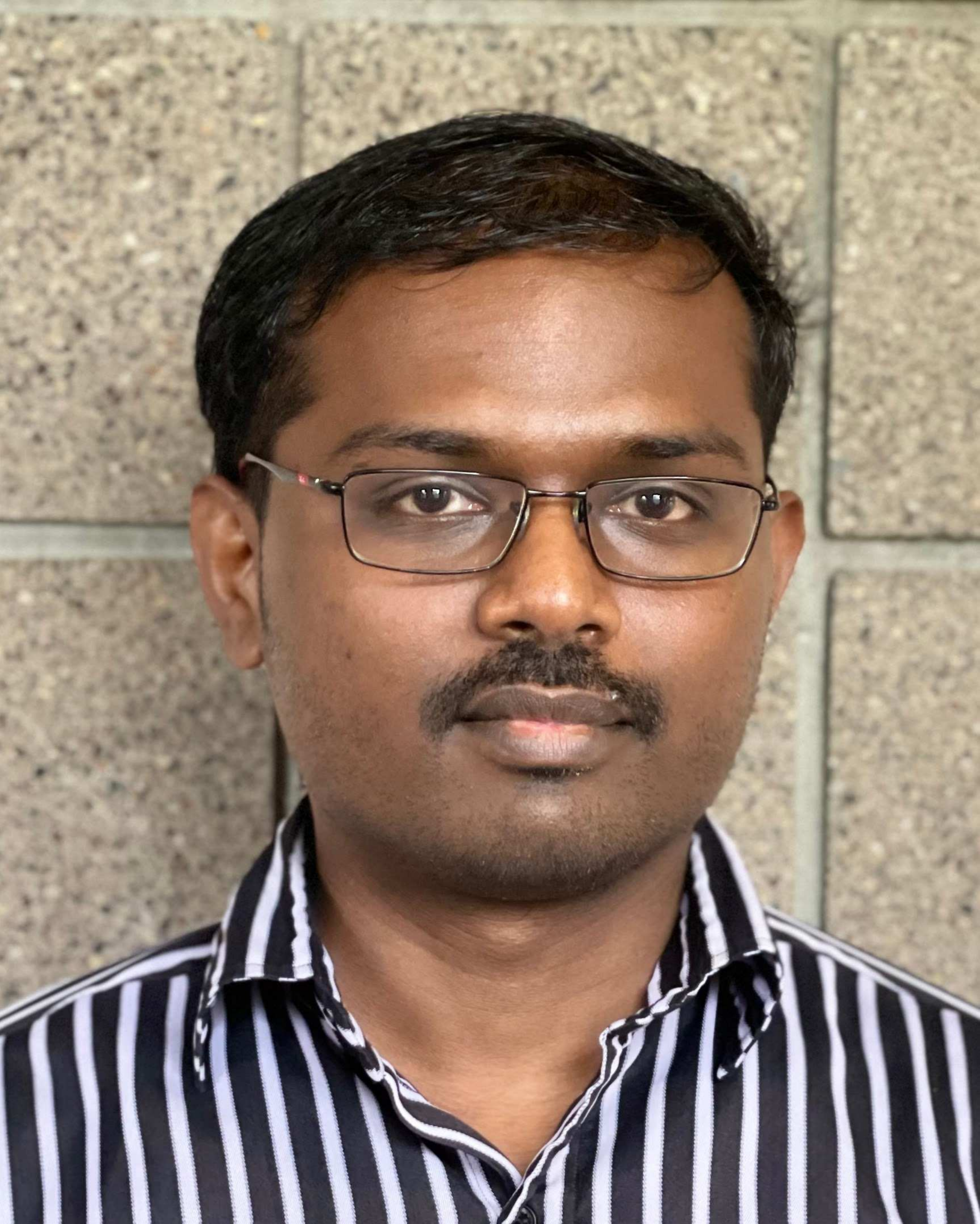}}]{Sethupathy Subramanian}
received the bachelors degree in electrical and electronics engineering from Anna University, Chennai, India, in 2009. He received the masters and doctrate degrees in electrical engineering from Indian Institute of Science, Bangalore, India in 2011 and 2017 respectively.
He is currently pursuing his graduate research at the Department of 
Physics, University of Notre Dame, USA.

His research interests, pertinent to electrical engineering, 
include computational electromagnetics, numerical stability, 
finite element and edge element methods. 

\end{IEEEbiography}
\begin{IEEEbiography}[{\includegraphics[width=1in,height=1.25in,clip,keepaspectratio]{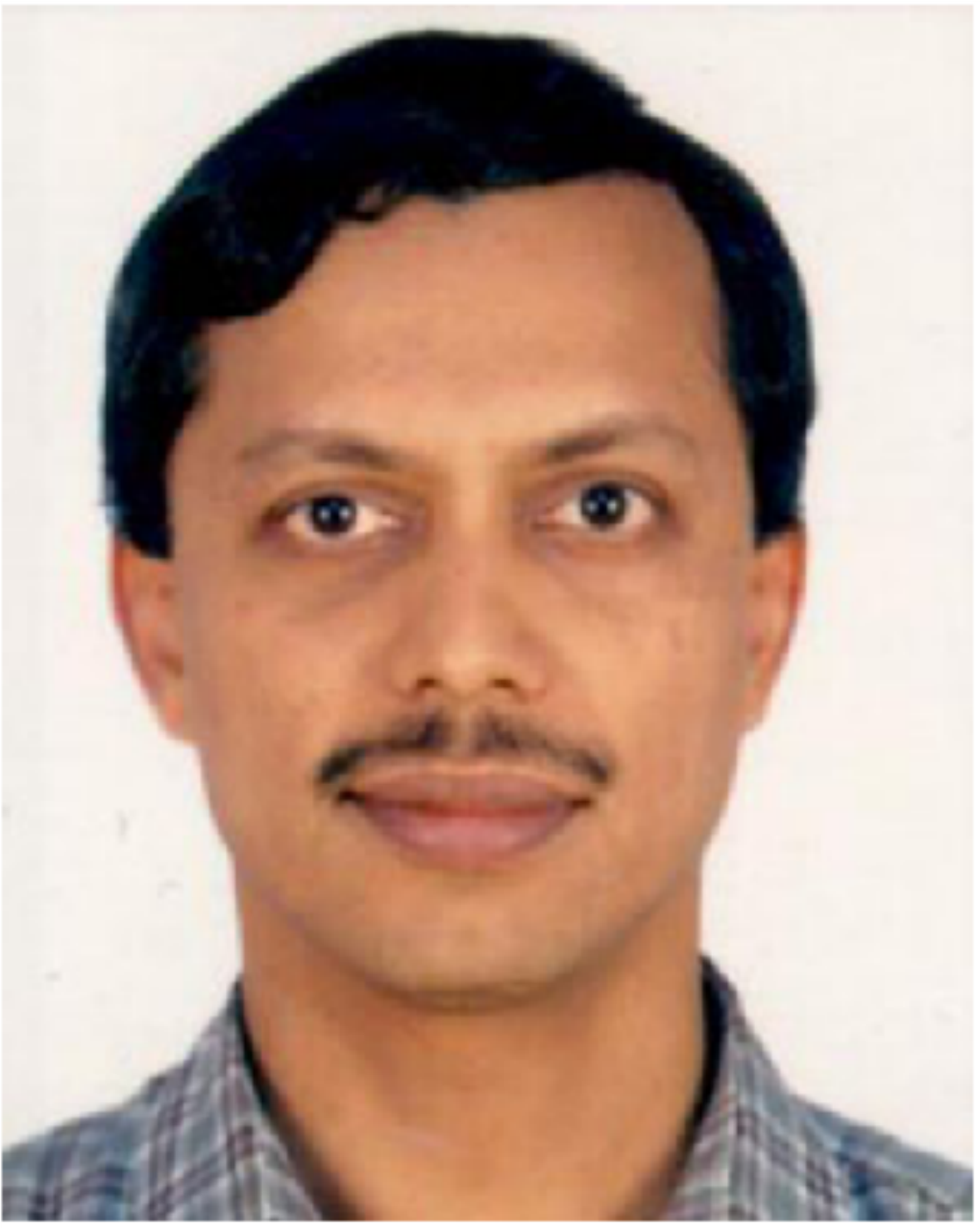}}]{Udaya Kumar}

 (Senior Member, IEEE) received the
bachelor’s degree in electrical engineering from 
Bangalore University, Bangaluru, India, in 1989, and the
M.E. and Ph.D. degrees in high-voltage engineering
from the Indian Institute of Science, Bengaluru,
India, in 1991 and 1998, respectively.

He is currently a Professor with the High Voltage
Laboratory, Department of Electrical Engineering,
Indian Institute of Science. He is an Associate Editor 
of IEEE Transactions on Power Delivery \&
IET High Voltage.

\end{IEEEbiography}

\begin{IEEEbiography}[{\includegraphics[width=1in,height=1.25in,clip,keepaspectratio]{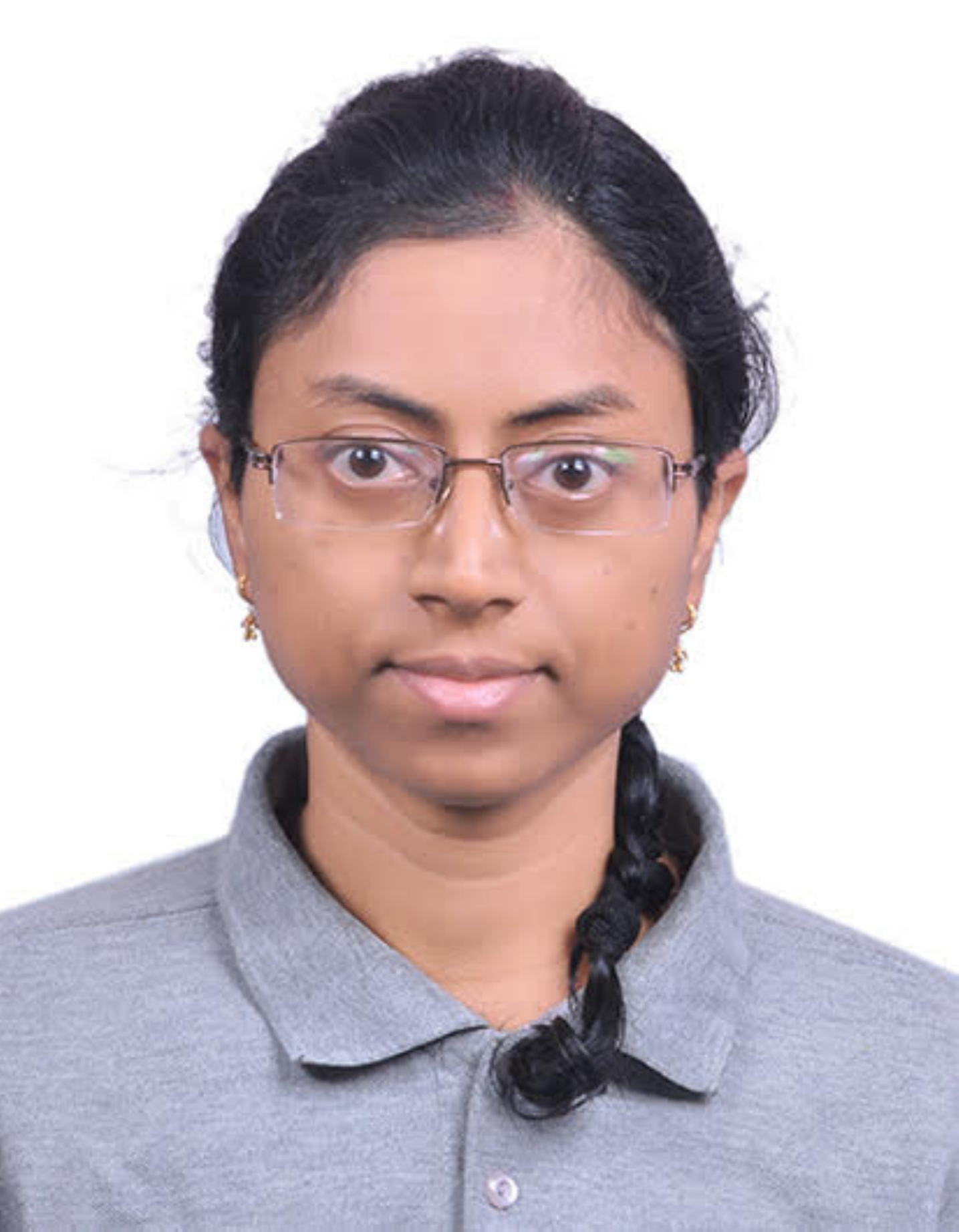}}]{Sujata Bhowmick}

received the B.E. degree in electrical engineering from IIEST, Shibpur, India, in 2006, and the M.E. degree in electrical engineering from the Indian Institute of Science, Bengaluru, India, in 2011. She received
the Ph.D. degree from the Department of Electronic Systems Engineering,
Indian Institute of Science, Bengaluru, India, in 2019.

Her current research interests include power electronics for renewable resources, single-phase grid-connected power converters, 
computational electromagnetics, finite element and 
edge element methods.

\end{IEEEbiography}

%
%




\end{document}